%% file: main.tex
\documentclass{article}

\usepackage{graphicx}
\usepackage{float}
\usepackage{amsthm}
\usepackage{amsmath}
\usepackage{amssymb}
\usepackage{mathrsfs}
\usepackage{tikz-cd}
\usepackage{stmaryrd}
\usepackage{enumitem}
\usepackage{placeins}
\usepackage{pdflscape}
\usepackage[colorlinks=true,citecolor=blue]{hyperref}
\usepackage{adjustbox}

\usepackage{marginnote}
\usepackage{comment}

\theoremstyle{definition}
\newtheorem{Definition}{Definition}[section]

\theoremstyle{definition}
\newtheorem{Notation}[Definition]{Notation}

\theoremstyle{plain}
\newtheorem{Theorem}[Definition]{Theorem}

\theoremstyle{plain}
\newtheorem{Proposition}[Definition]{Proposition}

\theoremstyle{definition}
\newtheorem{Construction}[Definition]{Construction}

\theoremstyle{plain}
\newtheorem{Lemma}[Definition]{Lemma}

\theoremstyle{plain}
\newtheorem{Corollary}[Definition]{Corollary}

\theoremstyle{plain}

\theoremstyle{definition}
\newtheorem{Example}[Definition]{Example}

\theoremstyle{remark}
\newtheorem{Remark}[Definition]{Remark}

\theoremstyle{plain}
\newcommand{\thistheoremname}{}
\newtheorem*{genericthm*}{\thistheoremname}
\newenvironment{namedthm*}[1]
  {\renewcommand{\thistheoremname}{#1}%
   \begin{genericthm*}}
  {\end{genericthm*}}

\title{Module Categories As Spans}
\author{Hao Xu\footnote{Email: \href{mailto:haoxu@imada.sdu.dk}{\tt haoxu@imada.sdu.dk}}}
\date{December 2025}

\begin{document}

\bibliographystyle{alpha}

\maketitle

\begin{abstract}
We establish a correspondence between modules and spans of algebras within a general monoidal 2-category $\mathfrak{C}$. Specifically, for an algebra $A$ in $\mathfrak{C}$, we construct a normalized lax 3-functor from the 2-category of $A$-modules to the 3-category of 2-spans of algebras in $\mathfrak{C}$ under $A$. This framework unifies and generalizes the realization of module functors and module natural transformations as spans of monoidal functors. We demonstrate the utility of this theory by recovering the realization of module objects in several familiar 2-categories and discuss its extension to the 2-categories $\mathbf{MCat}$ and $\mathbf{BrCat}$. In these cases, module objects correspond to central module monoidal categories over a braided monoidal category and central braided monoidal categories over a symmetric monoidal category, respectively.
\end{abstract}

{\hypersetup{linkcolor=black}\tableofcontents}

\section{Introduction}

Representation theory traditionally studies modules over associative algebras. A left module over an algebra $R$ can be defined in two equivalent ways:
\begin{itemize}
    \item A vector space $M$ together with a linear map $l:R \otimes M \to M$ satisfying associativity and unitality conditions: \[\begin{tikzcd}[column sep=large,row sep=large]
    {R \otimes R \otimes M} 
        \arrow[d,"\mathrm{id}_R \otimes l"'] 
        \arrow[r,"m^R \otimes \mathrm{id}_M"]
    & {R \otimes M}
        \arrow[d,"l"]
    \\ {R \otimes M}
        \arrow[r,"l"']
    & {M} \end{tikzcd}, \quad \begin{tikzcd}[column sep=huge,row sep=huge]
    {R \otimes M}
        \arrow[r, "l"]
    & {M}
        \arrow[d,equal]
    \\ {\Bbbk \otimes M}
        \arrow[u,"i^R \otimes \mathrm{id}_M"]
        \arrow[r,equal]
    & {M} \end{tikzcd};\]

    \item A vector space $M$ together with an algebra homomorphism $\phi:R \to \mathbf{End}(M)$, where $\mathbf{End}(M)$ is the algebra of endomorphisms on $M$.
\end{itemize} 
It is a standard observation that these definitions are interchangeable. Given the action $l$, one defines $\phi$ by $\phi(r)(m) = l(r \otimes m)$ for any $r \in R$ and $m \in M$. Conversely, given $\phi$, one recovers $l$.

However, modules over an algebra $R$ do not exist in isolation; they form a category where morphisms are $R$-linear maps. Two modules $M$ and $N$ are isomorphic if there exists an invertible $R$-linear map $f:M \to N$. In this case, $f$ induces an algebra isomorphism $f_*:\mathbf{End}(M) \simeq \mathbf{End}(N)$ compatible with the $R$-actions, i.e., $f_* \circ \phi^M = \phi^N$.

Characterizing non-invertible morphisms between modules in terms of the algebra actions is more subtle. A non-invertible linear map $f:M \to N$ does not induce an algebra homomorphism between endomorphism algebras. Instead, it induces two linear maps: \[f_* : \mathbf{End}(M) \to \mathbf{Hom}(M,N); \quad p \mapsto f \circ p,\] \[f^*: \mathbf{End}(N) \to \mathbf{Hom}(M,N); \quad q \mapsto q \circ f.\] 
The map $f$ preserves the $R$-actions if and only if the following diagram commutes: \[\begin{tikzcd}[column sep=large,row sep=large]
{R}
    \arrow[r,"\phi^M"]
    \arrow[d,"\phi^N"']
& {\mathbf{End}(M)}
    \arrow[d,"f_*"]
\\ {\mathbf{End}(N)}
    \arrow[r,"f^*"']  
& {\mathbf{Hom}(M,N)} \end{tikzcd}.\]

\noindent \textbf{Motivation.} The equivalence between these two definitions of a module, while elementary, highlights a fundamental dichotomy in how algebraic structures are conceptualized, which guides the methodology of this paper.

The first definition, which relies on the action map $l:R \otimes M \to M$, embodies the \emph{operadic perspective}. In this view, an algebra is essentially a monoid object in a monoidal category, and a module is an object equipped with an action by this monoid. This perspective naturally extends to algebras over operads, emphasizing internal composition laws and the geometrization of coherence data.

The second definition, based on the homomorphism $\phi:R \to \mathbf{End}(M)$, reflects the \emph{enriched perspective}. Here, the algebra $R$ is not merely an object but is regarded as a single-object category enriched over the base category (e.g., $\mathbf{Vect}$). A module, in this framework, is simply an enriched functor from this single-object category to the ambient enriched category.

While the operadic perspective is traditional and intuitive for defining modules, the enriched perspective offers a powerful computational advantage, especially in the context of higher category theory. By treating modules as functors, we can harness the robust machinery of \emph{enriched 2-category theory}. This shift in viewpoint not only simplifies computations but also provides a framework to externalize the structure of module functors, interpreting them as spans of algebras. Despite its potential, this enriched perspective remains underexplored in the literature, and this paper aims to bridge that gap by developing a comprehensive theory within the lower-dimensional setting of 2-categories. The insights gained here are expected to help interested readers explore these generalizations in the general framework of weak $(\infty,n)$-categories.

\subsection*{Main Results}

In this article, we work with a general semistrict monoidal 2-category $\mathfrak{C}$. We adopt the perspective that a module over an algebra $A$ in $\mathfrak{C}$ is an enriched 2-functor from the delooping $\mathbb{B}A$ into a chosen $\mathfrak{C}$-enriched 2-category $\mathfrak{M}$.

\renewcommand{\thistheoremname}{Definition \ref{def:ModulesInEnriched2Category}}

\begin{genericthm*}
We define the 2-category of \emph{$A$-modules in $\mathfrak{M}$} to be the 2-category of $\mathfrak{C}$-enriched 2-functors:  \[ \mathbf{Mod}_\mathfrak{M}(A) := \mathbf{Cat}(\mathfrak{C})(\mathbb{B} A, \mathfrak{M}). \]
\end{genericthm*}

To establish the notation, we develop the theory of pointed enriched 2-categories and the delooping construction in Section \ref{sec:Delooping}. The first main result is the functoriality of the delooping construction.

\renewcommand{\thistheoremname}{Proposition \ref{prop:DeloopingOfAlgebrasIs3Functor}}

\begin{genericthm*}
Delooping of algebras gives rise to a 3-functor \[ \mathbf{Alg}(\mathfrak{C}) \to {}^{\mathbf{I}/}_{oplax}\mathbf{Cat}(\mathfrak{C}), \] where $\mathbf{Alg}(\mathfrak{C})$ is viewed as a 3-category with only identity 3-morphisms. Moreover, this 3-functor induces equivalences on the level of hom 2-categories.
\end{genericthm*}

The second main result is the construction of a normalized lax 3-functor from the 2-category of modules into the 3-category of 2-spans of algebras.

\renewcommand{\thistheoremname}{Theorem \ref{thm:Embed2CatOfModulesInto2SpansOfAlgebras}}

\begin{genericthm*}
We construct a normalized lax 3-functor:
\[\mathbf{Mod}_{\mathfrak{M}}(A) \to \mathbf{Span}_{(3,2)}({}^{A/}\mathbf{Alg}(\mathfrak{C})).\]
\end{genericthm*}

This lax 3-functor provides a systematic dictionary between the theory of modules and the theory of spans:
\begin{itemize}
    \item An $A$-module $x$ is mapped to its enriched endomorphism algebra $\mathfrak{M}[x,x]$, regarded as an algebra under $A$ (Definition \ref{def:ModulesInEnriched2Category}).
    
    \item An $A$-module 1-morphism $f:x \to y$ is mapped to the span of algebras given by the following 2-fiber product (Theorem \ref{thm:AlgebraStructureOn2FiberProduct}): 
    \[\begin{tikzcd}
        {A_f}
            \arrow[d,"p_x"']
            \arrow[r,"p_y"]
        & {A_y}
            \arrow[d,"f^*"]
        \\ {A_x} 
            \arrow[r,"f_*"']
            \arrow[ur,Rightarrow,shorten <=10pt,shorten >=10pt,"\varpi"]
        & {\mathfrak{M}[x,y]}
    \end{tikzcd}.\]
    
    \item An $A$-module 2-morphism is mapped to a 2-span of algebras, constructed via comma objects (Theorem \ref{thm:2SpanOfAlgebrasAssociatedTo2Morphism}).
    
    \item The coherence data for the composition of 1-morphisms is provided in Theorem \ref{thm:Laxator1Morphism}.
    
    \item The coherence data for the composition of 2-morphisms is provided in Theorem \ref{thm:Laxator2Morphism}.
\end{itemize}

In the end of the article, we demonstrate the power of this general framework by applying it to concrete examples:

\begin{itemize}
    \item §\ref{sec:ApplicationToCat}: When $\mathfrak{C} = \mathbf{Cat}$ is the 2-category of small categories with the Cartesian product, we recover the functoriality that a module category over a monoidal category $\mathcal{C}$ corresponds to a monoidal functor from $\mathcal{C}$ into the endomorphism category of the module category.
    
    \item §\ref{sec:OtherExamples}: We list algebras and modules in several other examples of monoidal 2-categories, including:
    \begin{itemize}
        \item $\mathbf{2Vect}$, the 2-category of linear categories;
        
        \item $\mathbf{2Rep}(G)$, the 2-category of linear categories with an action of a finite group $G$;
        
        \item $\mathbf{2Vect}^\pi_G$, the 2-category of $G$-graded linear categories with a $\pi$-twisted monoidal structure, where $\pi \in \mathrm{H}^4(G,\Bbbk^\times)$ is a 4-cocycle on $G$;
        
        \item $\mathbf{Mod}(\mathcal{B})$, the 2-category of module categories over a braided fusion category $\mathcal{B}$.
    \end{itemize}
\end{itemize}

Beyond the scenarios covered by the assumptions of Theorem \ref{thm:Embed2CatOfModulesInto2SpansOfAlgebras}, we discuss in §\ref{sec:BeyondEnriched2Categories} how to extend the main theorem to the 2-categories $\mathbf{MCat}$ and $\mathbf{BrCat}$. In these cases, module objects correspond to central module monoidal categories over a braided monoidal category and central braided monoidal categories over a symmetric monoidal category, respectively. The former notion was studied in \cite{HPT16,HPT23} from the operadic perspective, while the latter notion has not appeared explicitly in the literature to the knowledge of the author.

\begin{Remark}
    The author was inspired by the work of Davydov, Kong and Runkel \cite{DKR} on formalizing bulk-boundary correspondence in two dimensional rational conformal field theories as cospans of algebras. It would be interesting to see if the ideas in this paper can be applied to study the bulk-boundary correspondence for general types of quantum field theories in higher dimensions.
\end{Remark}

\begin{Notation}
    Throughout this paper, we will use the following conventions:
    \begin{itemize}
        \item We follow the convention of Baez, where $n$-categories are always assumed to be weak, unless otherwise specified. Thus, bicategories (and other bicategorical structures such as bifunctors or pseudonatural transformations) in the literature will be denoted simply as 2-categories (2-functors, 2-natural transformations, respectively), while what people called 2-categories (and 2-functors) will be emphasized as strict 2-categories (and strict 2-functors).
        
        \item For 1-categorical notions, we often omit the prefix ``1-'' for brevity. 
        
        \item Since we are working with 2-categories, 3-categories and even 4-categories, it would be cumbersome to keep track of the levels of morphisms all the time. Therefore, we often use the single arrow $\to$ to denote morphisms of any level; the levels of morphisms will be clear from the source and target; in the pasting diagrams, we sometimes use double arrows $\Rightarrow$ to fill the 2-cells. We use $\mathbf{1}_x$ to denote the identity morphism on an object/morphism $x$, where the level of $\mathbf{1}_x$ is one level higher than that of $x$. 
        
        \item The composition of morphisms in the $k$-th direction is denoted by $\circ_k$. For examples, for 2-morphisms in a 2-category, $\circ_2$ denotes the vertical composition and $\circ_1$ denotes the horizontal composition. When there is no risk of confusion, we simply denote the composition by juxtaposition.
        
        \item Higher functors and natural transformations are by default strong, unless otherwise specified as lax or oplax.
        
        \item $\mathbf{Vect}$ denotes the category of finite dimensional vector spaces over the ground field $\Bbbk$, with the symmetric monoidal structure given by the tensor product of $\Bbbk$-vector spaces.
        
        \item For a monoidal category $\mathcal{C}$, we denote its monoidal product by $\otimes$, monoidal unit by $\mathbf{I}$, associator by $\alpha$, left unitor by $\lambda$, right unitor by $\rho$. When there is an additional braiding, we denote it by $\beta$. We add subscripts/superscripts to these symbols when we need to distinguish the structures from different monoidal categories.
        
        \item For a monoidal functor $F \colon \mathcal{C} \to \mathcal{D}$, we denote its multiplicative structure by $\gamma$; again we add subscripts/superscripts when necessary.
        
        \item We denote the \emph{monoidal opposite category} of a monoidal category $\mathcal{C}$ by $\mathcal{C}^{mp}$, and for a braided monoidal category $\mathcal{B}$, we denote the same monoidal category with the reversed braiding by $\mathcal{B}^{rev}$.
        
        \item While the definition of a fully weak monoidal 2-category involves a vast amount of coherence data, including 2-associators (pentagonators) and 2-unitors, we work primarily with \emph{semistrict monoidal 2-categories}, also known as \emph{Gray monoids}. A Gray monoid is defined as an algebra object in the symmetric monoidal category $\mathbf{Gray}$, consisting of \emph{strict} 2-categories and \emph{strict} 2-functors, and equipped with the \emph{Gray tensor product}. In a Gray monoid, the monoidal product is strictly associative and unital on objects, but the exchange law for 1-morphisms is relaxed to an invertible 2-morphism called the \emph{interchanger}. This approach is justified by the Coherence Theorem \cite{GPS}, which states that every monoidal 2-category is equivalent to a semistrict one. This allows us to suppress many coherence isomorphisms without loss of generality, significantly simplifying proofs.
        
        \item To manage the complexity of 2-categorical computations, we utilize the \emph{string diagram calculus}, which is rigorously justified by the strictification theorem for 2-categories. We adopt the convention following \cite{GS,D4} where: 
        \begin{itemize}
            \item Regions of the plane represent objects, and we omit labeling them for readability;
            
            \item Strings separating regions represent 1-morphisms, with composition depicted by placing strings in parallel from top to bottom;
            
            \item Nodes connecting strings represent 2-morphisms, with composition read from left to right.
        \end{itemize}

        For 1-morphisms $x \xrightarrow{f} y$, $u \xrightarrow{g} v$ in $\mathfrak{C}$, the \textit{interchanger} \[\begin{tikzcd}
            {x \, \Box \, u}
                \arrow[r,"x \, \Box \, g"]
                \arrow[d,"f \, \Box \, u"']
            & {x \, \Box \, v}
                \arrow[d,"f \, \Box \, v"]
                \arrow[dl,Rightarrow,shorten <= 10pt, shorten >= 10pt]
            \\ {y \, \Box \, u}
                \arrow[r,"y \, \Box \, g"']
            & {y \, \Box \, v}
        \end{tikzcd}\] is denoted by the string diagram below on the left,

        \newlength{\diagramlength}

        \settoheight{\diagramlength}{\includegraphics[width=20mm]{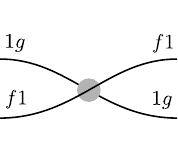}}

        \begin{center}
        \begin{tabular}{@{}ccc@{}}
            
        \includegraphics[width=20mm]{Pictures/Preliminaries/GraphicalCalculus/inter.pdf}  & \includegraphics[width=20mm]{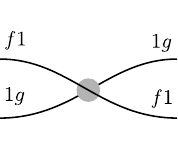}

        \end{tabular}
        \end{center} and we denote its inverse by the diagram on the right. Here $1$ is the place holder indicating which argument of the product $\Box$ is fixed. For readability, we abbreviate the notation $\Box$ in the string diagrams.
    \end{itemize}
\end{Notation}

\subsection*{Acknowledgements}

I would like to thank Thibault D{\'e}coppet for insightful discussions that led to the initial version of this article. I am also grateful to Liang Kong and Zhi-Hao Zhang for their valuable feedback on the early draft of this manuscript. This project was partially supported by the DAAD Graduate School Scholarship Programme (57572629) and the DFG Project 398436923.

This second version includes significant updates from the first. Sections 2, 5, and 6 of the first version have been reorganized into Section 5. The enriched 2-category framework, briefly mentioned in Section 5 of the first version, has been adopted as the central theme of this new version, following the suggestion of an anonymous referee. Reflecting this shift in narrative, Sections 2, 3, and 4 in the current version are new, where we carefully develop the arguments in the most general setting. I would like to thank Markus Zetto for his helpful comments on the draft of this second version.

\input{Preliminaries}

\input{Delooping}

\input{Span}

\input{Applications}

\bibliography{bibliography.bib}

\end{document}

%% file: Preliminaries.tex
\section{Preliminaries}

We fix a semistrict monoidal 2-category $\mathfrak{C}$.

\subsection{Algebras}

\begin{Definition} \label{def:Algebra}
    An \textit{algebra} in $\mathfrak{C}$ consists of:
    \begin{enumerate}
        \item An object $A$ in $\mathfrak{C}$;
        
        \item Two 1-morphisms, \textit{multiplication} $m: A \, \Box \, A \to A$ and \textit{unit} $i: \mathbf{I} \to A$;
        
        \item Three 2-isomorphisms, called \textit{associator} $\alpha$, \textit{left unitor} $\lambda$ and \textit{right unitor} $\rho$: \[\begin{tikzcd}
            {A \, \Box \, A \, \Box \, A} 
                \arrow[r, "m 1"]
                \arrow[d,"1 m"']
            & {A \, \Box \, A}
                \arrow[d, "m"]
                \arrow[ld,Rightarrow,shorten <= 10pt,shorten >= 10pt,"\alpha"]
            \\ {A \, \Box \, A}
                \arrow[r, "m"']
            & {A}
        \end{tikzcd}, \begin{tikzcd}
            {\mathbf{I} \, \Box \, A} 
                \arrow[r, "i 1"]
                \arrow[d,equal]
            & {A \, \Box \, A}
                \arrow[d, "m"]
                \arrow[dl,Rightarrow,shorten <= 10pt,shorten >= 10pt,"\lambda"']
                \arrow[dr,Rightarrow,shorten <= 10pt,shorten >= 10pt,"\rho"]
            & {A \, \Box \, \mathbf{I}} 
                \arrow[d,equal]
                \arrow[l,"1 i"']
            \\ {A}
                \arrow[r, equal]
            & {A}
            & {A}
                \arrow[l,equal]
        \end{tikzcd},\]
    \end{enumerate} subject to the following coherence conditions:
    
    \begin{enumerate}
        \item [a.] We have the pentagon equation in $\mathbf{Hom}_\mathfrak{C}(A \, \Box \, A \, \Box \, A \, \Box \, A,A)$:
    \end{enumerate}

    \settoheight{\diagramlength}{\includegraphics[height=24mm]{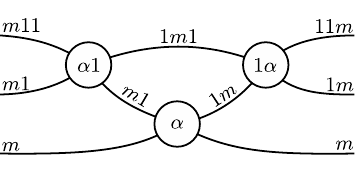}}

    \begin{equation} \label{eqn:AlgebraAssociativity}
        \begin{tabular}{@{}cccc@{}}
        
        \includegraphics[height=24mm]{Pictures/Preliminaries/Algebra/Algebra/Associativityleft.pdf} & \raisebox{0.45\diagramlength}{$=$} & \includegraphics[height=24mm]{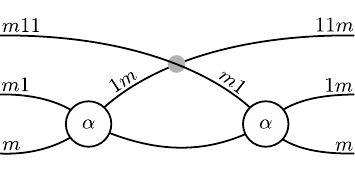} & \raisebox{0.45\diagramlength}{,}
        
        \end{tabular}
    \end{equation}

    \begin{enumerate}
        \item [b.] We have the triangle equation in $\mathbf{Hom}_\mathfrak{C}(A \, \Box \, A,A)$:
    \end{enumerate}

    \settoheight{\diagramlength}{\includegraphics[height=24mm]{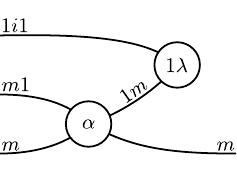}}

    \begin{equation} \label{eqn:AlgebraUnitality}
        \begin{tabular}{@{}cccc@{}}

        \includegraphics[height=24mm]{Pictures/Preliminaries/Algebra/Algebra/Unitalityleft.pdf} & \raisebox{0.45\diagramlength}{$=$} &
        \includegraphics[height=24mm]{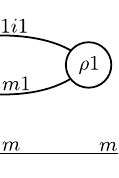} & \raisebox{0.45\diagramlength}{.}

        \end{tabular}
    \end{equation}
\end{Definition}

\subsection{Algebra 1-Morphisms}

\begin{Definition} \label{def:Algebra1Morphism}
    Let $A$ and $B$ be two algebras in $\mathfrak{C}$. An \textit{algebra 1-morphism} $f:A\rightarrow B$ consists of a 1-morphism $f:A\rightarrow B$ in $\mathfrak{C}$, together with two 2-isomorphisms
     
    \begin{center}
        \begin{tabular}{@{}c c@{}}
            \begin{tikzcd}[sep=small]
                {A \, \Box \, A} 
                    \arrow[rr, "f1"] 
                    \arrow[dd, "m^A"'] 
                &  {}
                & {B \, \Box \, A} 
                    \arrow[rr, "1f"] 
                & {} 
                & {B \, \Box \, B} 
                    \arrow[dd, "m^B"] 
                    \arrow[lllldd, Rightarrow, "\psi^f"', shorten > = 15mm, shorten < = 15mm] 
                \\ {}
                & {}
                & {}
                & {} 
                & {}
                \\ {A} 
                    \arrow[rrrr, "f"']
                & {}
                & {} 
                & {} 
                & {B}
            \end{tikzcd},
            &
            \begin{tikzcd}[sep=small]
                {} 
                & {} 
                & {A}
                    \arrow[rrdd, "f"]
                & {} 
                & {}  
                \\ {}
                & {} 
                & {}  
                \\ {}
                {\mathbf{I}} 
                    \arrow[rruu, "i^A"] 
                    \arrow[rrrr, "i^B"'] 
                & {} 
                & {} 
                    \arrow[uu, Rightarrow, "\eta^f"', shorten > = 2ex, shorten < = 1ex] 
                & {} 
                & {B}
            \end{tikzcd},
        \end{tabular}
    \end{center}
    satisfying:
    
    \begin{enumerate}
    \item [a.] We have
    \end{enumerate}
    
    \begin{center}
    \begin{equation} \label{eqn:Algebra1MorphismCoh1}
    \begin{tabular}{@{}ccc@{}}
    
    \includegraphics[height=40mm]{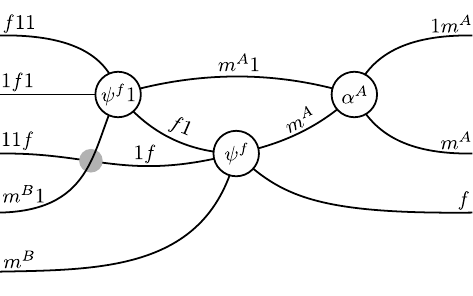} \\ 
    \rotatebox{90}{$=$} \\
    \includegraphics[height=40mm]{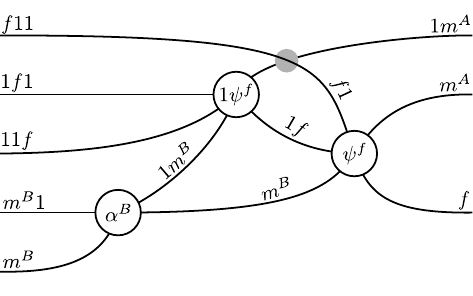}
    
    \end{tabular}
    \end{equation}
    \end{center}
    
    \begin{enumerate}
        \item [] in $\mathbf{Hom}_\mathfrak{C}(A \, \Box \, A \, \Box \, A,B)$;
    
        \item [b.] We have
    \end{enumerate}
    
    \settoheight{\diagramlength}{\includegraphics[height=24mm]{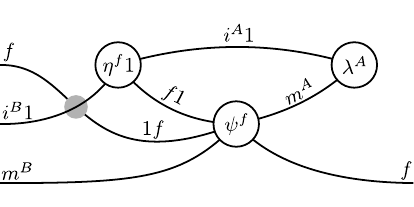}}
    
    \begin{center}
    \begin{equation} \label{eqn:Algebra1MorphismCoh2}
    \begin{tabular}{@{}ccc@{}}
    
    \includegraphics[height=28mm]{Pictures/Preliminaries/Algebra/Algebra1Morphism/algebra1morcoh2left.pdf} & \raisebox{0.45\diagramlength}{$=$} &
    
    \includegraphics[height=24mm]{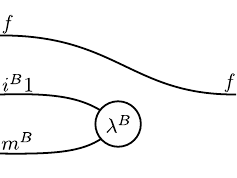}
    
    \end{tabular}
    \end{equation}
    \end{center}
    
    \begin{enumerate}
        \item [] in $\mathbf{Hom}_\mathfrak{C}(A,B)$;
        
        \item [c.] We have
    \end{enumerate}
    
    \settoheight{\diagramlength}{\includegraphics[height=24mm]{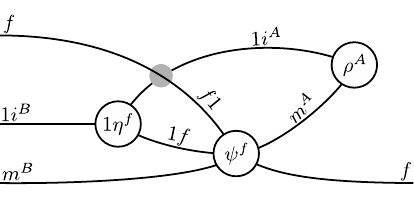}}
    
    \begin{center}
    \begin{equation} \label{eqn:Algebra1MorphismCoh3}
    \begin{tabular}{@{}ccc@{}}
    
    \includegraphics[height=28mm]{Pictures/Preliminaries/Algebra/Algebra1Morphism/algebra1morcoh3left.pdf} & \raisebox{0.45\diagramlength}{$=$} &
    
    \includegraphics[height=24mm]{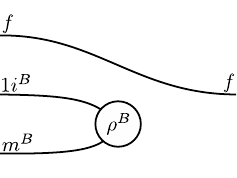}
    
    \end{tabular}
    \end{equation}
    \end{center}

    \begin{enumerate}
        \item [] in $\mathbf{Hom}_\mathfrak{C}(A,B)$.
    \end{enumerate}
\end{Definition}

\subsection{Algebra 2-Morphisms}
    
\begin{Definition} \label{def:Algebra2Morphism}
    Let $A$ and $B$ be two algebras, $f$ and $g$ be two algebra 1-morphisms from $A$ to $B$ in $\mathfrak{C}$. An \textit{algebra 2-morphism} $\gamma$ consists of a 2-morphism $\gamma:f \to g$ in $\mathfrak{C}$, satisfying:
    \begin{enumerate}
    \item [a.] We have
    \end{enumerate}
    
    \settoheight{\diagramlength}{\includegraphics[height=24mm]{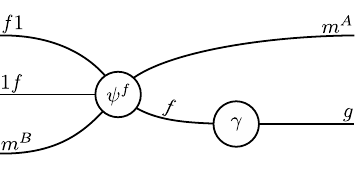}}
    
    \begin{center}
    \begin{equation} \label{eqn:Algebra2MorphismCoh1}
    \begin{tabular}{@{}ccc@{}}
    
    \includegraphics[height=24mm]{Pictures/Preliminaries/Algebra/Algebra2Morphism/algebra2morcoh1left.pdf} & \raisebox{0.45\diagramlength}{$=$} &
    
    \includegraphics[height=24mm]{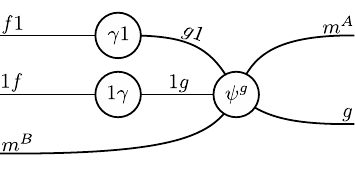}
    
    \end{tabular}
    \end{equation}
    \end{center}
    
    \begin{enumerate}
        \item [] in $\mathbf{Hom}_\mathfrak{C}(A \, \Box \, A,B)$;
    
        \item [b.] We have
    \end{enumerate}
    
    \settoheight{\diagramlength}{\includegraphics[height=24mm]{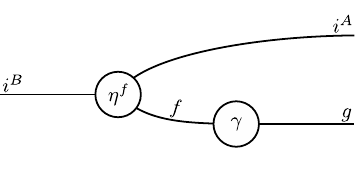}}
    
    \begin{center}
    \begin{equation} \label{eqn:Algebra2MorphismCoh2}
    \begin{tabular}{@{}ccc@{}}
    
    \includegraphics[height=24mm]{Pictures/Preliminaries/Algebra/Algebra2Morphism/algebra2morcoh2left.pdf} & \raisebox{0.45\diagramlength}{$=$} &
    
    \includegraphics[height=24mm]{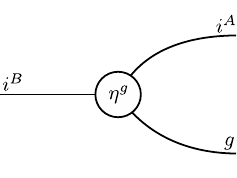}
    
    \end{tabular}
    \end{equation}
    \end{center}

    \begin{enumerate}
        \item [] in $\mathbf{Hom}_\mathfrak{C}(\mathbf{I},B)$.
    \end{enumerate}
\end{Definition}
    
\begin{Construction}
    There is a 2-category $\mathbf{Alg}(\mathfrak{C})$ where:
    \begin{itemize}
        \item Objects are algebras in $\mathfrak{C}$;
        
        \item 1-morphisms are algebra 1-morphisms in $\mathfrak{C}$;
        
        \item 2-morphisms are algebra 2-morphisms in $\mathfrak{C}$.
    \end{itemize}
\end{Construction}

\subsection{Enriched 2-Categories}

In this section, we recall the basic notions of enriched 2-categories from \cite{GS}. The concept of enriched categories was first introduced by B{\'e}nabou \cite{Ben63} and Kelly \cite{Kel82} independently to generalize ordinary categories by replacing the Hom-sets with objects in a monoidal category. This idea was later extended to 2-categories in two theses \cite{Car95,Lac95} in the same year.

\begin{Definition} \label{def:Enriched2Category}
    An \textit{enriched 2-category} over a semistrict monoidal 2-category $\mathfrak{C}$ consists of:
    \begin{enumerate}
        \item A set\footnote{We assume smallness in general and ignore the size issue.} of objects $Ob(\mathfrak{M})$;

        \item For objects $x,y$ in $Ob(\mathfrak{M})$, an \textit{enriched Hom}\footnote{When there is no ambiguity, we shall omit the $\mathfrak{M}$ and $\mathfrak{C}$ from the notation.} object $\mathfrak{M}[x,y]_\mathfrak{C}$ in $\mathfrak{C}$;

        \item For object $x$ in $Ob(\mathfrak{M})$, an \textit{enriched unit} $j_x:\mathbf{I} \to \mathfrak{M}[x,x]_\mathfrak{C}$;

        \item For objects $x,y,z$ in $Ob(\mathfrak{M})$, an \textit{enriched composition} \[m_{x,y,z}:\mathfrak{M}[y,z]_\mathfrak{C} \,\Box \, \mathfrak{M}[x,y]_\mathfrak{C} \to \mathfrak{M}[x,z]_\mathfrak{C};\]

        \item For objects $x,y$ in $Ob(\mathfrak{M})$, \textit{enriched unitors} in $\mathfrak{C}$ \[\begin{tikzcd}[sep=30pt]
            {\mathbf{I} \, \Box \, \mathfrak{M}[x,y]_\mathfrak{C}}
                \arrow[r,equal]
                \arrow[d,"j_y 1"']
            & {\mathbf{I} \, \Box \, \mathfrak{M}[x,y]_\mathfrak{C}}
                \arrow[d,equal]
            \\ {\mathfrak{M}[y,y]_\mathfrak{C} \, \Box \, \mathfrak{M}[x,y]_\mathfrak{C}}
                \arrow[r,"m_{x,y,y}"']
                \arrow[ur,Rightarrow,shorten <=20pt, shorten >=20pt,"\lambda_{x,y}"]
            & {\mathfrak{M}[x,y]_\mathfrak{C}}
        \end{tikzcd},\] \[\begin{tikzcd}[sep=30pt]
            {\mathfrak{M}[x,y]_\mathfrak{C} \, \Box \, \mathbf{I}}
                \arrow[r,equal]
                \arrow[d,"1 j_x"']
            & {\mathfrak{M}[x,y]_\mathfrak{C} \, \Box \, \mathbf{I}}
                \arrow[d,equal]
            \\ {\mathfrak{M}[x,y]_\mathfrak{C} \, \Box \, \mathfrak{M}[x,x]_\mathfrak{C}}
                \arrow[r,"m_{x,x,y}"']
                \arrow[ur,Rightarrow,shorten <=20pt, shorten >=20pt,"\rho_{x,y}"]
            & {\mathfrak{M}[x,y]_\mathfrak{C}}
        \end{tikzcd};\]

        \item For objects $x,y,z,w$ in $Ob(\mathfrak{M})$, an \textit{enriched associator} in $\mathfrak{C}$ \[\begin{tikzcd}[sep=30pt]
            {(\mathfrak{M}[z,w]_\mathfrak{C} \, \Box \, \mathfrak{M}[y,z]_\mathfrak{C}) \, \Box \, \mathfrak{M}[x,y]_\mathfrak{C}}
                \arrow[r,"m_{y,z,w} 1"]
                \arrow[d,equal]
            & {\mathfrak{M}[y,w]_\mathfrak{C} \, \Box \, \mathfrak{M}[x,y]_\mathfrak{C}}
                \arrow[dd,"m_{x,y,w}"]
            \\{\mathfrak{M}[z,w]_\mathfrak{C} \, \Box \,( \mathfrak{M}[y,z]_\mathfrak{C} \, \Box \, \mathfrak{M}[x,y]_\mathfrak{C})}
                \arrow[d,"1 m_{x,y,z}"']
            & {}
                \arrow[l,Rightarrow,shorten <=10pt, shorten >=10pt,"\pi_{x,y,z,w}"']
            \\ {\mathfrak{M}[z,w]_\mathfrak{C} \, \Box \, \mathfrak{M}[x,z]_\mathfrak{C}}
                \arrow[r,"m_{x,z,w}"']
            & {\mathfrak{M}[x,w]_\mathfrak{C}}
        \end{tikzcd};\]
    \end{enumerate}
    subject to the following simplified conditions in \cite[Section 3.1]{GS}:
    \begin{enumerate}
        \item [a.] For objects $x,y,z,u,v$ in $Ob(\mathfrak{M})$, we have
    \end{enumerate}
    
    \settoheight{\diagramlength}{\includegraphics[height=24mm]{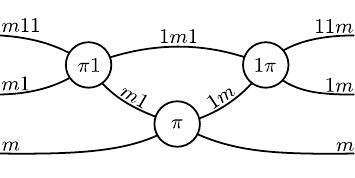}}
    
    \begin{equation}\label{eqn:Enriched2CategoryAssociativity}
    \begin{tabular}{@{}ccc@{}}
    
    \includegraphics[height=24mm]{Pictures/Preliminaries/Enriched2Category/EnrichedAssociativityleft.pdf} & \raisebox{0.45\diagramlength}{$=$} &
    \includegraphics[height=24mm]{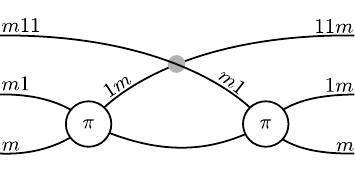}
    
    \end{tabular}
    \end{equation}

    \begin{enumerate}
        \item [] in $\mathbf{Hom}_{\mathfrak{C}}([u,v] \, \Box \, [z,u] \, \Box \, [y,z] \, \Box \, [x,y], [x,v])$,
        \item [b.] For objects $x,y,z$ in $Ob(\mathfrak{M})$, we have
    \end{enumerate}

    \settoheight{\diagramlength}{\includegraphics[height=24mm]{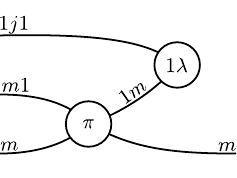}}
    
    \begin{equation}\label{eqn:Enriched2CategoryUnitality}
    \begin{tabular}{@{}ccc@{}}
    
    \includegraphics[height=24mm]{Pictures/Preliminaries/Enriched2Category/EnrichedUnitalityleft.pdf} & \raisebox{0.45\diagramlength}{$=$} &
    \includegraphics[height=24mm]{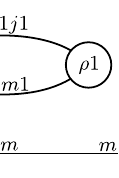}
    
    \end{tabular}
    \end{equation}

    \begin{enumerate}
        \item [] in $\mathbf{Hom}_{\mathfrak{C}}([y,z] \, \Box \, [x,y], [x,z])$.
    \end{enumerate}
\end{Definition}
    
\begin{Construction} \label{cstr:Underlying2CategoryOfEnriched2Category}
    By post-composing the enriched Hom objects with the forgetful 2-functor
    \[
    \mathfrak{C} \to \mathbf{Cat}; \quad x \mapsto \mathbf{Hom}_\mathfrak{C}(\mathbf{I},x),
    \]
    we obtain a $\mathbf{Cat}$-enriched 2-category structure on $\mathfrak{M}$. This construction yields the \textit{underlying 2-category} of $\mathfrak{M}$, whose objects are precisely $Ob(\mathfrak{M})$. By slight abuse of notation, we continue to denote the underlying 2-category by $\mathfrak{M}$.

    A 1-morphism $f \colon x \to y$ in the underlying 2-category $\mathfrak{M}$ corresponds to a 1-morphism $f^\sharp \colon \mathbf{I} \to  \mathfrak{M}[x,y]_\mathfrak{C}$ in $\mathfrak{C}$. Similarly, a 2-morphism $a \colon f \to g$ in the underlying 2-category $\mathfrak{M}$ corresponds to a 2-morphism $a^\sharp \colon f^\sharp \to g^\sharp$ in $\mathfrak{C}$.
\end{Construction}

\subsection{Enriched Hom-2-Functors}

\begin{Construction} \label{cstr:EnrichedHom2Functor}
    Enriched Hom can be promoted to a strict cubical 2-functor \[\mathfrak{M}^{1op} \times \mathfrak{M} \to \mathfrak{C}; \quad (x,y) \mapsto \mathfrak{M}[x,y]_\mathfrak{C}.\]
    
    \noindent More precisely, for any object $x$ in $\mathfrak{M}$, there is a 2-functor \[\mathfrak{M}[x,-]: \mathfrak{M} \to \mathfrak{C}\] sending \begin{itemize}
        \item An object $y$ in $\mathfrak{M}$ to the object $\mathfrak{M}[x,y]$ in $\mathfrak{C}$;
        
        \item A 1-morphism $f:y \to z$ in $\mathfrak{M}$ to the composite \[\mathfrak{M}[x,y] = \mathbf{I} \, \Box \, \mathfrak{M}[x,y] \xrightarrow{f^\sharp 1} \mathfrak{M}[y,z] \, \Box \, \mathfrak{M}[x,y] \xrightarrow{m^\mathfrak{M}_{x,y,z}} \mathfrak{M}[x,z], \] which we denote by $\mathfrak{M}[1,f]:\mathfrak{M}[x,y] \to \mathfrak{M}[x,z]$ or simply $f_*$;
        
        \item A 2-morphism $a:f \to f'$, where $f,f':y \to z$ are 1-morphisms in $\mathfrak{M}$, to the composite \[\begin{tikzcd}
            {\mathfrak{M}[x,y]}
                \arrow[r,equal]
            & {\mathbf{I} \, \Box \, \mathfrak{M}[x,y]}
                \arrow[r,"f^\sharp 1" {name=U},bend left=20pt]
                \arrow[r,"{f^\prime}^\sharp 1"' {name=D},bend right=20pt]
                \arrow[from=U,to=D,Rightarrow,shorten <=5pt, shorten >=5pt,"a^\sharp 1"']
            & {\mathfrak{M}[y,z] \, \Box \, \mathfrak{M}[x,y]}
                \arrow[r,"m^\mathfrak{M}_{x,y,z}"]
            & {\mathfrak{M}[x,z]}
        \end{tikzcd},\] which we denote by $\mathfrak{M}[1,a]:\mathfrak{M}[1,f] \to \mathfrak{M}[1,f']$ or simply $a_*$. 
    \end{itemize}

    \noindent Similarly, for any object $y$ in the $\mathfrak{M}$, there is a 2-functor \[\mathfrak{M}[-,y]: \mathfrak{M}^{1op} \to \mathfrak{C}\] sending \begin{itemize}
        \item An object $x$ in $\mathfrak{M}$ to the object $\mathfrak{M}[x,y]$ in $\mathfrak{C}$;
        
        \item A 1-morphism $g:w \to x$ in $\mathfrak{M}$ to the composite \[\mathfrak{M}[x,y] = \mathfrak{M}[x,y] \, \Box \, \mathbf{I} \xrightarrow{1 g^\sharp} \mathfrak{M}[x,y] \, \Box \, \mathfrak{M}[w,x] \xrightarrow{m^\mathfrak{M}_{w,x,y}} \mathfrak{M}[w,y], \] which we denote by $\mathfrak{M}[g,1]:\mathfrak{M}[x,y] \to \mathfrak{M}[w,y]$ or simply $g^*$;
        
        \item A 2-morphism $c:g \to g'$, where $g,g':w \to x$ are 1-morphisms in $\mathfrak{M}$, to the composite \[\begin{tikzcd}
            {\mathfrak{M}[x,y]}
                \arrow[r,equal]
            & {\mathfrak{M}[x,y] \, \Box \, \mathbf{I}}
                \arrow[r,"1 g^\sharp" {name=U},bend left=20pt]
                \arrow[r,"1 {g^\prime}^\sharp"' {name=D},bend right=20pt]
                \arrow[from=U,to=D,Rightarrow,shorten <=5pt, shorten >=5pt,"1 c^\sharp"']
            & {\mathfrak{M}[x,y] \, \Box \, \mathfrak{M}[w,x]}
                \arrow[r,"m^\mathfrak{M}_{w,x,y}"]
            & {\mathfrak{M}[w,y]}
        \end{tikzcd},\] which we denote by $\mathfrak{M}[c,1]:\mathfrak{M}[g,1] \to \mathfrak{M}[g',1]$ or simply $c^*$.
    \end{itemize}

    \noindent Moreover, for 1-morphisms $y \xrightarrow{f} z$ and $w \xrightarrow{g} x$ in $\mathfrak{M}$, there is an interchanger of enriched Hom 2-functors \[ \mathfrak{M}[g,1] \circ \mathfrak{M}[1,f] \to \mathfrak{M}[1,f] \circ \mathfrak{M}[g,1] \] induced from the interchanger $\pmb{\phi}_{f,g}$ in $\mathfrak{C}$.

    Finally, enriched associator $\pi^\mathfrak{M}$ and enriched unitors $\lambda^\mathfrak{M},\rho^\mathfrak{M}$ are promoted to be 2-natural isomorphisms between the corresponding enriched Hom 2-functors.
\end{Construction}

\subsection{Enriched 2-Functors}
Let $\mathfrak{M}$ and $\mathfrak{N}$ be two $\mathfrak{C}$-enriched 2-categories.

\begin{Definition} \label{def:Enriched2Functor}
    A \textit{$\mathfrak{C}$-enriched 2-functor} from $\mathfrak{M}$ to $\mathfrak{N}$ consists of:
    \begin{enumerate}
        \item A function $F:Ob(\mathfrak{M}) \to Ob(\mathfrak{N})$;
        
        \item For objects $x,y$ in $Ob(\mathfrak{M})$, a 1-morphism $F_{x,y}: \mathfrak{M}[x,y] \to \mathfrak{N}[Fx,Fy]$;
        
        \item For objects $x,y,z$ in $Ob(\mathfrak{M})$, a 2-isomorphism $m^F_{x,y,z}$ in $\mathfrak{C}$ \[\begin{tikzcd}[column sep=40pt]
            {\mathfrak{M}[y,z] \, \Box \, \mathfrak{M}[x,y]}
                \arrow[d,"F_{y,z} \, \Box \, F_{x,y}"']
                \arrow[r,"m^\mathfrak{M}_{x,y,z}"]
            & {\mathfrak{M}[x,z]}
                \arrow[d,"F_{x,z}"]
            \\ {\mathfrak{N}[Fy,Fz] \, \Box \, \mathfrak{N}[Fx,Fy]}
                \arrow[r,"m^\mathfrak{N}_{Fx,Fy,Fz}"']
                \arrow[ur,Rightarrow,shorten <=20pt, shorten >=20pt,"m^F_{x,y,z}"]
            & {\mathfrak{N}[Fx,Fz]}
        \end{tikzcd};\]

        \item For object $x$ in $Ob(\mathfrak{M})$, a 2-isomorphism $j^F_x$ in $\mathfrak{C}$ \[\begin{tikzcd}
            {\mathbf{I}}
                \arrow[r,"j^\mathfrak{M}_x"]
                \arrow[d,equal]
            & {\mathfrak{M}[x,x]}
                \arrow[d,"F_{x,x}"]
            \\ {\mathbf{I}}
                \arrow[r,"j^\mathfrak{N}_{Fx}"']
                \arrow[ur,Rightarrow,shorten <=10pt, shorten >=10pt,"j^F_x"]
            & {\mathfrak{N}[Fx,Fx]}
        \end{tikzcd};\]
    \end{enumerate} subject to the following conditions:

    \begin{enumerate}
        \item[a.] For objects $x,y,z,w$ in $Ob(\mathfrak{M})$, we have
    \end{enumerate}

    \settoheight{\diagramlength}{\includegraphics[height=40mm]{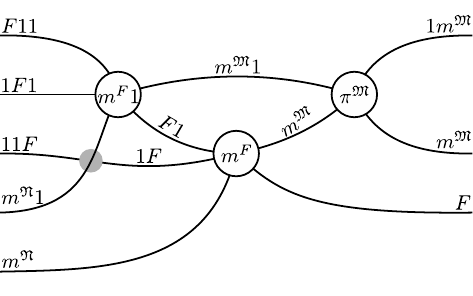}}

    \begin{equation}\label{eqn:Enriched2FunctorAssociativity}
    \begin{tabular}{@{}c@{}}
        \includegraphics[height=40mm]{Pictures/Preliminaries/Enriched2Functor/Enriched2FunctorAssociativityLeft.pdf} \\ {\rotatebox{90}{$=$}} \\
        \includegraphics[height=40mm]{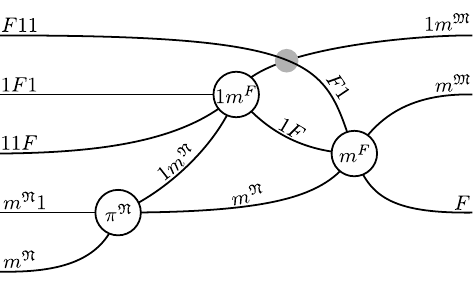}
        \end{tabular}
    \end{equation}

    \begin{enumerate}
        \item [] in $\mathbf{Hom}_{\mathfrak{C}}(\mathfrak{M} [z,w] \, \Box \, \mathfrak{M} [y,z] \, \Box \, \mathfrak{M} [x,y], \mathfrak{N}[Fx,Fw])$,
        \item[b.] For objects $x,y$ in $Ob(\mathfrak{M})$, we have
    \end{enumerate}

    \settoheight{\diagramlength}{\includegraphics[height=24mm]{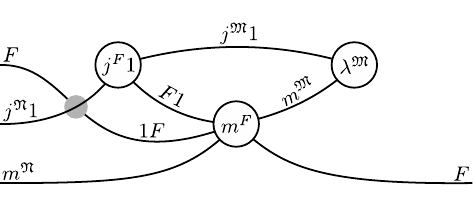}}

    \begin{equation}\label{eqn:Enriched2FunctorUnitality1}
    \begin{tabular}{@{}ccc@{}}
        \includegraphics[height=24mm]{Pictures/Preliminaries/Enriched2Functor/Enriched2FunctorUnitality1Left.pdf} & \raisebox{0.45\diagramlength}{$=$} &
        \includegraphics[height=24mm]{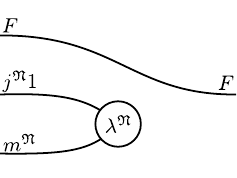}
        \end{tabular}
    \end{equation}

    \begin{enumerate}
        \item [] in $\mathbf{Hom}_{\mathfrak{C}}(\mathfrak{M}[x,y], \mathfrak{N}[Fx,Fy])$,
        \item[c.] For objects $x,y$ in $Ob(\mathfrak{M})$, we have
    \end{enumerate}

    \settoheight{\diagramlength}{\includegraphics[height=24mm]{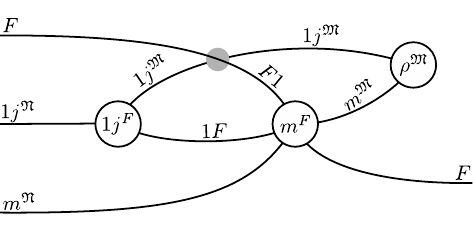}}

    \begin{equation}\label{eqn:Enriched2FunctorUnitality2}
    \begin{tabular}{@{}ccc@{}}
        \includegraphics[height=24mm]{Pictures/Preliminaries/Enriched2Functor/Enriched2FunctorUnitality2Left.pdf} & \raisebox{0.45\diagramlength}{$=$} &
        \includegraphics[height=24mm]{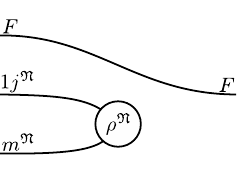}
        \end{tabular}
    \end{equation}

    \begin{enumerate}
        \item [] in $\mathbf{Hom}_{\mathfrak{C}}(\mathfrak{M}[x,y], \mathfrak{N}[Fx,Fy])$.
    \end{enumerate}
\end{Definition}

\begin{Construction} \label{cstr:Underlying2FunctorOfEnriched2Functor}
    Every $\mathfrak{C}$-enriched 2-functor $F:\mathfrak{M} \to \mathfrak{N}$ induces a 2-functor between the underlying 2-categories via the following:
    \begin{itemize}
        \item The underlying 2-functor $F:\mathfrak{M} \to \mathfrak{N}$ is given by the same function on objects;
        
        \item For objects $x,y$ in $\mathfrak{M}$, we have a functor induced between Hom categories \[\mathbf{Hom}_\mathfrak{C}(\mathbf{I},F_{x,y}): \mathbf{Hom}_\mathfrak{C}(\mathbf{I},\mathfrak{M}[x,y]) \to \mathbf{Hom}_\mathfrak{C}(\mathbf{I},\mathfrak{N}[Fx,Fy]). \]
    \end{itemize}

    For simplicity, for 1-morphism $f$ and 2-morphism $a$ in $\mathfrak{M}$, we still denote their image under the above functor as $Ff$ and $Fa$, respectively. The 2-functor $F$ is called the \textit{underlying 2-functor} of the enriched 2-functor.
\end{Construction}

\begin{Remark}
    Enriched 2-functors and enriched Hom 2-functors are compatible. More precisely, there are canonically induced coherence data: \[\begin{tikzcd}
        {\mathfrak{M}[x,y]}
            \arrow[d,"{\mathfrak{M}[1,f]}"']
            \arrow[r,"F_{x,y}"]
        & {\mathfrak{N}[Fx,Fy]}
            \arrow[d,"{\mathfrak{N}[1,Ff]}"]
            \arrow[ld,Rightarrow,shorten <=10pt, shorten >=10pt,"F_{x,f}"]
        \\ {\mathfrak{M}[x,z]}
            \arrow[r,"F_{x,z}"']
        & {\mathfrak{N}[Fx,Fz]}
    \end{tikzcd}, \quad \begin{tikzcd}
        {\mathfrak{M}[x,y]}
            \arrow[d,"{\mathfrak{M}[g,1]}"']
            \arrow[r,"F_{x,y}"]
        & {\mathfrak{N}[Fx,Fy]}
            \arrow[d,"{\mathfrak{N}[Fg,1]}"]
            \arrow[ld,Rightarrow,shorten <=10pt, shorten >=10pt,"F_{g,y}"]
        \\ {\mathfrak{M}[w,y]}
            \arrow[r,"F_{w,y}"']
        & {\mathfrak{N}[Fw,Fy]}
    \end{tikzcd},\] for 1-morphisms $y \xrightarrow{f} z$ and $w \xrightarrow{g} x$ in $\mathfrak{M}$.
\end{Remark}

\subsection{Enriched 2-Natural Transformations}

Let $F,G:\mathfrak{M} \to \mathfrak{N}$ be two $\mathfrak{C}$-enriched 2-functors.

\begin{Definition} \label{def:Enriched2NatTrans}
    A \textit{$\mathfrak{C}$-enriched 2-natural transformation} from $F$ to $G$ consists of:
    \begin{enumerate}
        \item For object $x$ in $\mathfrak{M}$, a 1-morphism $\eta^\sharp_x:\mathbf{I} \to \mathfrak{N}[Fx,Gx]$ in $\mathfrak{C}$, or equivalently, a 1-morphism $\eta_x: Fx \to Gx$ in the underlying 2-category of $\mathfrak{N}$;
        
        \item For objects $x,y$ in $\mathfrak{M}$, a 2-isomorphism in $\mathfrak{C}$: \[\begin{tikzcd}[column sep=40pt]
            {\mathfrak{M}[x,y]}
                \arrow[r,"F_{x,y}"]
                \arrow[d,"G_{x,y}"']
            & {\mathfrak{N}[Fx,Fy]}
                \arrow[d,"{\mathfrak{N}[1,\eta_y]}"]
            \\ {\mathfrak{N}[Gx,Gy]}
                \arrow[r,"{\mathfrak{N}[\eta_x,1]}"']
                \arrow[ur,Rightarrow,shorten <=15pt, shorten >=15pt,"\eta_{x,y}"]
            & {\mathfrak{N}[Fx,Gy]}
        \end{tikzcd};\]
    \end{enumerate} subject to the following conditions:

    \begin{enumerate}
        \item [a.] For objects $x,y,z$ in $\mathfrak{M}$, we have
    \end{enumerate}

    \settoheight{\diagramlength}{\includegraphics[height=32mm]{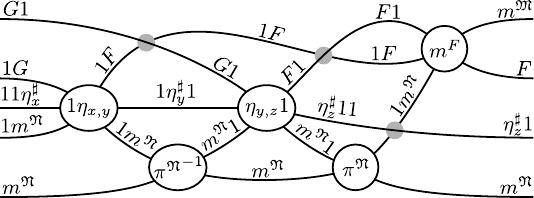}}

    \begin{equation} \label{eqn:Enriched2NatTransAssociativity}
    \begin{tabular}{@{}c@{}}
        \includegraphics[height=32mm]{Pictures/Preliminaries/Enriched2NatTrans/AssociativeLeft.pdf} \\ {\rotatebox{90}{$=$}} \\
        \includegraphics[height=32mm]{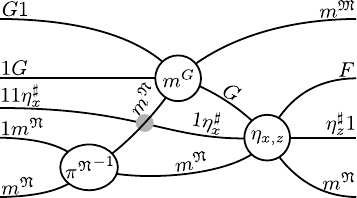}
        \end{tabular}
    \end{equation}

    \begin{enumerate}
        \item [] in $\mathbf{Hom}_{\mathfrak{C}}(\mathfrak{M}[y,z] \, \Box \, \mathfrak{M}[x,y], \mathfrak{N}[Fx,Gz])$;
        
        \item [b.] For object $x$ in $\mathfrak{M}$, we have
    \end{enumerate}

    \settoheight{\diagramlength}{\includegraphics[height=24mm]{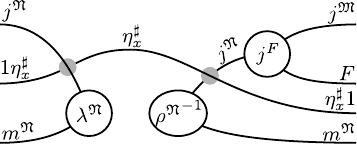}}

    \begin{equation} \label{eqn:Enriched2NatTransUnitality}
    \begin{tabular}{@{}ccc@{}}
        \includegraphics[height=24mm]{Pictures/Preliminaries/Enriched2NatTrans/UnitalLeft.pdf} & \raisebox{0.45\diagramlength}{$=$} &
        \includegraphics[height=24mm]{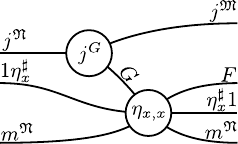}
        \end{tabular}
    \end{equation}

    \begin{enumerate}
        \item [] in $\mathbf{Hom}_{\mathfrak{C}}(\mathbf{I}, \mathfrak{N}[Fx,Gx])$.        
    \end{enumerate}

\end{Definition}

\begin{Construction} \label{cstr:Underlying2NatTrans}
    Every $\mathfrak{C}$-enriched 2-natural transformation $\eta:F \to G$ induces a 2-natural transformation between the underlying 2-functors:
    \begin{itemize}
        \item For object $x$ in $\mathfrak{M}$, $\eta$ assigns the same 1-morphism $\eta_x:Fx \to Gx$ in the underlying 2-category of $\mathfrak{N}$;
        
        \item For 1-morphism $x \xrightarrow{f} y$ in $\mathfrak{M}$, $\eta$ assigns the 2-isomorphism $\eta_f := \eta_{x,y} \circ_1 \mathbf{1}_{f^\sharp}$ in $\mathbf{Hom}_\mathfrak{C}(\mathbf{I},\mathfrak{N}[Fx,Gy])$.
    \end{itemize}
\end{Construction}

\begin{Remark}
    In Definition \ref{def:Enriched2NatTrans}, if we do not require $\eta_{x,y}$ to be invertible, then we obtain the notion of \textit{lax $\mathfrak{C}$-enriched 2-natural transformation}. Similarly, if we reverse the direction of $\eta_{x,y}$, then we obtain the notion of \textit{oplax $\mathfrak{C}$-enriched 2-natural transformation}.
\end{Remark}

\subsection{Enriched Modifications}

Let $\eta,\xi:F \to G$ be two $\mathfrak{C}$-enriched 2-natural transformations between $\mathfrak{M}$ and $\mathfrak{N}$.

\begin{Definition} \label{def:EnrichedModification}
    A \textit{$\mathfrak{C}$-enriched modification} from $\eta$ to $\xi$ consists of: a 2-morphism $\theta_x: \eta_x \to \xi_x$ in $\mathfrak{N}$ for each object $x$ in $\mathfrak{M}$, such that for objects $x,y$ in $\mathfrak{C}$, the equation 
    
    \settoheight{\diagramlength}{\includegraphics[height=24mm]{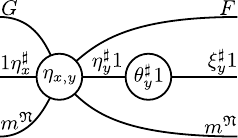}}

    \begin{equation} \label{eqn:EnrichedModification}
    \begin{tabular}{@{}ccc@{}}
        \includegraphics[height=24mm]{Pictures/Preliminaries/EnrichedModification/Left.pdf} & \raisebox{0.45\diagramlength}{$=$} &
        \includegraphics[height=24mm]{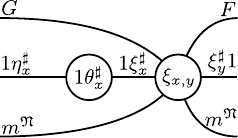}
        \end{tabular}
    \end{equation}
    
    \noindent holds in $\mathbf{Hom}_\mathfrak{C}(\mathfrak{M}[x,y],\mathfrak{N}[Fx,Gy])$.
\end{Definition}

\begin{Construction} \label{cstr:UnderlyingModification}
    Every $\mathfrak{C}$-enriched modification $\theta$ induces a modification between the underlying 2-natural transformations, where each component $\theta_x$ is the same 2-morphism for object $x$ in $\mathfrak{M}$.
\end{Construction}

\begin{Construction}
    For a given monoidal 2-category $\mathfrak{C}$, we obtain a 3-category $\mathbf{Cat}(\mathfrak{C})$ consists of:
    \begin{itemize}
        \item $\mathfrak{C}$-enriched 2-categories as objects,
        
        \item $\mathfrak{C}$-enriched 2-functors as 1-morphisms,
        
        \item $\mathfrak{C}$-enriched 2-natural transformations as 2-morphisms,
        
        \item $\mathfrak{C}$-enriched modifications as 3-morphisms.
    \end{itemize}

    There is a canonical forgetful 3-functor \[\mathbf{Cat}(\mathfrak{C}) \to \mathbf{2Cat}\] forgetting the $\mathfrak{C}$-enrichment, as in Construction \ref{cstr:Underlying2CategoryOfEnriched2Category}, \ref{cstr:Underlying2FunctorOfEnriched2Functor}, \ref{cstr:Underlying2NatTrans} and \ref{cstr:UnderlyingModification}.
\end{Construction}

\begin{Example}
    Ordinary 2-categories are just $\mathbf{Cat}$-enriched 2-categories.
\end{Example}

\subsection{2-Fiber Product and Comma Object}

\begin{Definition} \label{def:2FiberProduct}
Given a cospan in the 2-category $\mathfrak{C}$, its \emph{2-fiber product}
\[\begin{tikzcd}[sep=large]
    {X \times_Z Y}
        \arrow[d,"\mathbf{pr2}"']
        \arrow[r,"\mathbf{pr1}"]
    & {X}
        \arrow[d,"f"]
        \arrow[dl,Rightarrow,shorten <=25pt, shorten >=20pt,"\phi"']
    \\ {Y}
        \arrow[r,"g"']
    & {Z}
\end{tikzcd}\] has the following 2-universal property:
\begin{enumerate}
    \item For any object $A$ in $\mathfrak{C}$, 1-morphisms $p \colon A \to X$, $q \colon A \to Y$ and 2-isomorphism $\xi:f \circ p \simeq g \circ q$, there exists a 1-morphism $u \colon A \to X \times_Z Y$ and 2-isomorphisms $\zeta_1 \colon \mathbf{pr_1} \circ u \simeq p$ and $\zeta_2 \colon \mathbf{pr_2} \circ u \simeq q$, such that \[\xi \circ_2 (\mathbf{1}_f \circ_1 \zeta_1) = (\mathbf{1}_g \circ_1 \zeta_2) \circ_2 (\phi \circ_1 \mathbf{1}_u);\]

    \item For any 1-morphisms $u,v \colon A \to X \times_Z Y$ and 2-morphisms 
    \[\gamma_1 \colon \mathbf{pr_1} \circ u \to \mathbf{pr_1} \circ v, \quad \gamma_2 \colon \mathbf{pr_2} \circ u \to \mathbf{pr_2} \circ v,\] 
    such that 
    \[(\phi \circ_1 \mathbf{1}_v) \circ_2 (\mathbf{1}_f \circ_1 \gamma_1) = (\mathbf{1}_g \circ_1 \gamma_2) \circ_2 (\phi \circ_1 \mathbf{1}_u),\]
    then there exists a unique 2-morphism $\theta \colon u \to v$ such that $\gamma_1 = \mathbf{1}_\mathbf{pr_1} \circ_1 \theta$ and $\gamma_2 = \mathbf{1}_\mathbf{pr_2} \circ_1 \theta$.
\end{enumerate}

\begin{Construction} \label{cstr:CommaObject}
    For our purpose, we need a weaker version of 2-fiber product, where we do not require $\phi$ (and $\xi$ respectively) to be invertible. We denote this using the notation $X {}_f;_g Y$, which is inspired by the classical comma category construction.
\end{Construction}

\end{Definition}

%% file: Delooping.tex
\section{Delooping} \label{sec:Delooping}

\subsection{Pointed enriched 2-categories}

\begin{Definition}
    A \emph{pointed $\mathfrak{C}$-enriched 2-category} is a $\mathfrak{C}$-enriched 2-category $\mathfrak{M}$ equipped with a $\mathfrak{C}$-enriched 2-functor $x: \mathbf{I} \to \mathfrak{M}$. Equivalently, $\mathfrak{M}$ is equipped with a distinguished object $x$ in its underlying 2-category.
\end{Definition}

\begin{Definition}
    We define a 3-category ${}^{\mathbf{I}/}_{oplax}\mathbf{Cat}(\mathfrak{C})$ where:
    \begin{itemize}
        \item An object is a pointed $\mathfrak{C}$-enriched 2-category;
        
        \item A 1-morphism from $(\mathfrak{M},x)$ to $(\mathfrak{N},y)$ is a \emph{pointed $\mathfrak{C}$-enriched 2-functor}, which consists of a $\mathfrak{C}$-enriched 2-functor $F \colon \mathfrak{M} \to \mathfrak{N}$ and a $\mathfrak{C}$-enriched 2-natural isomorphism $\alpha \colon F \circ x \simeq y$;

        \item A 2-morphism from $(F,\alpha)$ to $(G,\beta)$ is a \emph{pointed $\mathfrak{C}$-enriched oplax 2-natural transformation}, which consists of a $\mathfrak{C}$-enriched oplax 2-natural transformation $\eta: F \to G$ and an invertible $\mathfrak{C}$-enriched modification $\mu: \beta \circ_2 (\eta \circ_1 \mathbf{1}_x) \to \alpha$;

        \item A 3-morphism from $(\eta,\mu)$ to $(\xi,\nu)$ is a \emph{pointed $\mathfrak{C}$-enriched modification}, which consists of a $\mathfrak{C}$-enriched modification $\theta: \eta \to \xi$ such that \[\nu \circ_3 (\beta \circ_2 (\theta \circ_1 \mathbf{1}_x)) = \mu.\]
    \end{itemize}
    In particular, if we consider pointed enriched \emph{strong} 2-natural transformations instead of icons, this 3-category becomes the (co-)slice 3-category of objects in $\mathbf{Cat}(\mathfrak{C})$ under $\mathbf{I}$.
\end{Definition}

\begin{Remark}
    Lack's notion of an \emph{icon} \cite{Lac10} is an oplax 2-natural transformation between 2-functors that is identity on objects. Up to equivalence of 2-categories, an icon can be identified with an oplax 2-natural transformation whose components at objects are invertible. When 2-categories are regarded as enriched over $\mathbf{Cat}$, our notion of pointed oplax 2-natural transformation is weaker than Lack's icon: we require only the component at the distinguished object to be invertible. In the special case where the target has a single object, the two notions agree.
\end{Remark}

\subsection{Delooping of Algebras}

\begin{Construction}
    For each algebra $A$ in $\mathfrak{C}$, we can construct its \emph{delooping} $\mathbb{B} A$ as a $\mathfrak{C}$-enriched 2-category with one single object, such that its enriched endo-hom recovers the algebra $A$. Comparing Definition \ref{def:Algebra} and Definition \ref{def:Enriched2Category}, we see that $\mathbb{B} A$ is equipped with a canonical pointed $\mathfrak{C}$-enriched 2-category structure. 
\end{Construction}

\begin{Proposition} \label{prop:DeloopingOfAlgebrasIs3Functor}
    Delooping of algebras gives rise to a 3-functor \[ \mathbf{Alg}(\mathfrak{C}) \to {}^{\mathbf{I}/}_{oplax}\mathbf{Cat}(\mathfrak{C}), \] where $\mathbf{Alg}(\mathfrak{C})$ is viewed as a 3-category with only identity 3-morphisms. Moreover, this 3-functor induces equivalences on the level of hom 2-categories.
\end{Proposition}

\begin{proof}
    On the level of objects, this is immediate: Definition \ref{def:Algebra} is the special case of Definition \ref{def:Enriched2Category} with a single object. Similarly, on the level of 1-morphisms, Definition \ref{def:Algebra1Morphism} coincides with Definition \ref{def:Enriched2Functor}, and the pointing is preserved automatically.
    
    For algebras $A,B$ in $\mathfrak{C}$, we need a 2-functor
    \[
        \mathbf{Alg}(\mathfrak{C})(A,B) \to {}^{\mathbf{I}/}_{oplax}\mathbf{Cat}(\mathfrak{C})(\mathbb{B} A,\mathbb{B} B).
    \]
    On 1-morphisms, Definition \ref{def:Algebra2Morphism} agrees with Definition \ref{def:Enriched2NatTrans}. The extra requirement that the component of the enriched oplax 2-natural transformation at the unique object be isomorphic to \(\mathbf{1}_B\) removes the redundancy.
    
    On 2-morphisms, the source has only trivial 2-morphisms since \(\mathbf{Alg}(\mathfrak{C})\) admits only identity 2-morphisms. In the target, a pointed enriched modification between deloopings of algebras reduces to a single condition because the pointing determines the relevant data. Therefore the constructed 3-functor induces equivalences on hom 2-categories.
\end{proof}

\subsection{Equivalent Characterization of Modules}

Suppose $\mathfrak{M}$ is a $\mathfrak{C}$-enriched 2-category and $A$ is an algebra in $\mathfrak{C}$.

\begin{Definition} \label{def:ModulesInEnriched2Category}
    We define the 2-category of \emph{$A$-modules in $\mathfrak{M}$} to be the 2-category of $\mathfrak{C}$-enriched 2-functors:  \[ \mathbf{Mod}_\mathfrak{M}(A) := \mathbf{Cat}(\mathfrak{C})(\mathbb{B} A, \mathfrak{M}). \]
\end{Definition}

\begin{Remark}
    Unpacking the definitions, an $A$-module in $\mathfrak{M}$ consists of an object $x$ in $\mathfrak{M}$, which is the image of the unique object under the enriched 2-functor, together with an algebra 1-morphism $l^x: A \to \mathfrak{M}[x,x]$ in $\mathfrak{C}$.

    Let $(x,l^x)$ and $(y,l^y)$ be two $A$-modules in $\mathfrak{M}$. An $A$-module 1-morphism from $(x,l^x)$ to $(y,l^y)$ consists of a 1-morphism $f \colon x \to y$ in $\mathfrak{M}$ together with a 2-isomorphism in $\mathfrak{C}$: \[\begin{tikzcd}
        {A}
            \arrow[r,"l^y"]
            \arrow[d,"l^x"']
        & {\mathfrak{M}[y,y]}
            \arrow[d,"f^*"]
        \\ {\mathfrak{M}[x,x]}
            \arrow[r,"f_*"']
            \arrow[ur,Rightarrow,shorten <=10pt, shorten >=10pt,"l^f"]  
        & {\mathfrak{M}[x,y]}
    \end{tikzcd}\] subject to three extra coherence conditions, as instances of \eqref{eqn:Enriched2FunctorAssociativity}, \eqref{eqn:Enriched2FunctorUnitality1} and \eqref{eqn:Enriched2FunctorUnitality2} where there is a single object.

    Let $(f,l^f)$ and $(g,l^g)$ be two $A$-module 1-morphisms between $(x,l^x)$ and $(y,l^y)$. An $A$-module 2-morphism consists of a 2-morphism $a \colon f \to g$ in $\mathfrak{M}$, subject to two coherence conditions coming from \eqref{eqn:Enriched2NatTransAssociativity} and \eqref{eqn:Enriched2NatTransUnitality}.
\end{Remark}

\begin{Remark}
    When $\mathfrak{C}$ is a closed monoidal 2-category, i.e. a monoidal 2-category with internal homs satisfying the expected 2-adjunction, the notion of an $A$-module in $\mathfrak{M} = \mathfrak{C}$ recovers the notion of a left module over a pseudomonoid as introduced in \cite{JS91}. In particular, when $\mathfrak{C} = \mathbf{Cat}$, we recover the classical notion of a module category over a monoidal category introduced by Bernstein \cite{Ber95}.
\end{Remark}

%% file: Span.tex
\section{Presenting Enriched 2-Categories As 2-Spans} \label{sec:ConstructionOfLax3Functor}

Suppose $\mathfrak{M}$ is a $\mathfrak{C}$-enriched 2-category where $\mathfrak{C}$ is a monoidal 2-category with 2-fiber products and comma objects. In this section, we will show that there is a canonical lax 3-functor from $\mathfrak{M}$ into $\mathbf{Span}_{(3,2)}(\mathbf{Alg}(\mathfrak{C}))$.

\subsection{Spans of Algebras}

\begin{Definition}
    There is a 4-category of 2-spans of algebras in $\mathfrak{C}$, denoted as $\mathbf{Span}_{(4,2)}(\mathbf{Alg}(\mathfrak{C}))$, defined as follows:
    \begin{itemize}
        \item Objects are algebras in $\mathfrak{C}$;
        
        \item 1-morphisms are spans of algebras in $\mathfrak{C}$;
        
        \item 2-morphisms are 2-spans between spans of algebras in $\mathfrak{C}$, e.g., for algebras $B$ and $C$ in $\mathfrak{C}$, a 2-span between $B$ and $C$ is a span of algebras over $B \times C$, or more explicitly, a diagram of the form \[\begin{tikzcd}
            {}
            & {}
            & {D}
                \arrow[lldd,"f_0"']
                \arrow[rrdd,"g_0"]
            & {}
            & {}
            \\ {}
            & {}
                \arrow[dd,Rightarrow,shorten <=15pt,shorten >=15pt,"\phi^f"]
            & {}
            & {}
                \arrow[dd,Rightarrow,shorten <=15pt,shorten >=15pt,"\phi^g"']
            & {}
            \\ {B}
            & {}
            & {E}
                \arrow[uu,"p" description]
                \arrow[dd,"q" description]
            & {}
            & {C}
            \\ {}
            & {}
            & {}
            & {}
            & {}
            \\ {}
            & {}
            & {F}
                \arrow[lluu,"f_1"]
                \arrow[rruu,"g_1"']
            & {}
            & {}
        \end{tikzcd}\] where $f_0, f_1, g_0, g_1, p, q$ are algebra 1-morphisms, $\phi^f, \phi^g$ are algebra 2-isomorphisms;

        \item 3-morphisms are 1-morphisms between 2-spans, e.g., for algebras $B$ and $C$ in $\mathfrak{C}$, 1-spans of algebras $(D,f_0,G_0)$ and $(F,f_1,g_1)$ between $B$ and $C$, a 1-morphism between 2-spans $(E_0,p_0,q_0,\phi^f_0,\phi^g_0)$ and $(E_1,p_1,q_1,\phi^f_1,\phi^g_1)$ consists of the following data:
        \[\begin{tikzcd}
            {}
            & {}
                \arrow[dd,Rightarrow,shorten <=20pt,shorten >=20pt,"\phi^p"]
            & {E_0}
                \arrow[lld,"p_0"']
                \arrow[rrd,"q_0"]
                \arrow[dd,"h" description]
            & {}
                \arrow[dd,Rightarrow,shorten <=20pt,shorten >=20pt,"\phi^q"']
            & {}
            \\ {D}
            & {}
            & {E}
            & {}
            & {F}
            \\ {}
            & {}
            & {E_1}
                \arrow[llu,"p_1"]
                \arrow[rru,"q_1"']
            & {}
            & {}
        \end{tikzcd}\] where $h$ is an algebra 1-morphism, $\phi^p, \phi^q$ are algebra 2-isomorphisms, subject to the obvious compatibility conditions;

        \item 4-morphisms are 2-morphisms between 2-spans, e.g., given 3-morphisms $(h_0,\phi^p_0,\phi^q_0)$ and $(h_1,\phi^p_1,\phi^q_1)$ as above, a 4-morphism between them is an algebra 2-morphism from $h_0$ to $h_1$ in $\mathfrak{C}$ such that the obvious compatibility conditions hold;
        
        \item Composition of 1-morphisms is given by the 2-fiber products \[\begin{tikzcd}
            {}
            & {}
            & {F \times_D G}
                \arrow[ld,"g^!(k)"']
                \arrow[rd,"k^!(g)"]
            & {}
            & {}
            \\ {}
            & {F}
                \arrow[ld,"f"']
                \arrow[rd,"g"']
                \arrow[rr,Rightarrow,shorten <=23pt,shorten >=20pt]
            & {}
            & {G}
                \arrow[ld,"k"]
                \arrow[rd,"l"]
            & {}
            \\ {C}
            & {}
            & {D}
            & {}
            & {E}
        \end{tikzcd};\]
        
        \item Composition of 2-morphisms is also induced by the 2-fiber product of the two vertical spans in the given composable 2-spans;
        
        \item Composition of 3-morphisms and 4-morphisms are induced by the compositions of the underlying algebra 1-morphisms and 2-morphisms in $\mathfrak{C}$ respectively.
    \end{itemize}
\end{Definition}

\begin{Construction}
    We can truncate the 4-category $\mathbf{Span}_{(4,2)}(\mathbf{Alg}(\mathfrak{C}))$ to a 3-category $\mathbf{Span}_{(3,2)}(\mathbf{Alg}(\mathfrak{C}))$ by identifying 3-morphisms related by invertible 4-morphisms. Similarly, one can further truncate $\mathbf{Span}_{(3,2)}(\mathbf{Alg}(\mathfrak{C}))$ to a 2-category $\mathbf{Span}_{(2,2)}(\mathbf{Alg}(\mathfrak{C}))$ by identifying 2-morphisms related by invertible 3-morphisms.
\end{Construction}

\subsection{On the level of Objects}

We assign each object $x$ in $\mathfrak{M}$ to its enriched endo-hom object $A_x := \mathfrak{M}[x,x]$, which is an algebra in $\mathfrak{C}$.
    
\subsection{On the level of 1-Morphisms}

Recall from Construction \ref{cstr:EnrichedHom2Functor} that for each 1-morphism $f \colon x \to y$ in $\mathfrak{M}$, there is associated push-forward $f_*$ and pullback $f^*$ between enriched hom objects. Consider the 2-fiber product \[\begin{tikzcd}
    {A_f}
        \arrow[d,"p_x"']
        \arrow[r,"p_y"]
    & {A_y}
        \arrow[d,"f^*"]
    \\ {A_x} 
        \arrow[r,"f_*"']
        \arrow[ur,Rightarrow,shorten <=10pt,shorten >=10pt,"\varpi"]
    & {\mathfrak{M}[x,y]}
\end{tikzcd}.\]

\begin{Theorem} \label{thm:AlgebraStructureOn2FiberProduct}
    There is a canonical algebra structure on $A_f$ such that \[A_x \xleftarrow{p_x} A_f \xrightarrow{p_y} A_y\] can be promoted to a span of algebras in $\mathfrak{C}$.
\end{Theorem}

\begin{proof}
    First, let us provide the multiplication on $A_f$, which should be a 1-morphism $m_f \colon A_f \Box A_f \to A_f$ in $\mathfrak{C}$. By the 2-universal property of 2-fiber product, it suffices to provide the following data:
    \begin{itemize}
        \item 1-morphisms \[A_f \Box A_f \xrightarrow{p_x \Box p_x} A_x \Box A_x \xrightarrow{m_{x,x,x}} A_x\] and \[A_f \Box A_f \xrightarrow{p_y \Box p_y} A_y \Box A_y \xrightarrow{m_{y,y,y}} A_y;\]
        
        \item one 2-isomorphism given by the following composite
    \end{itemize}

    \settoheight{\diagramlength}{\includegraphics[height=40mm]{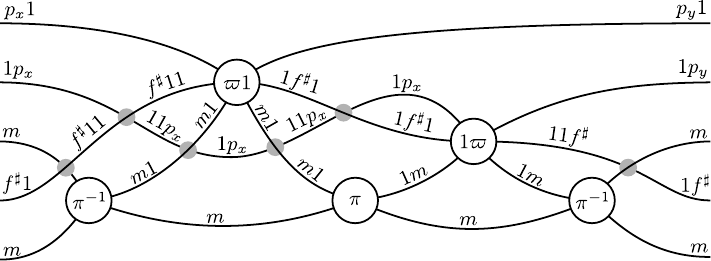}}

    \begin{equation*} \label{eqn:1MorphismMultiplication}
    \begin{tabular}{@{}ccc@{}}
        \includegraphics[height=40mm]{Pictures/1Morphism/multiplication.pdf} \raisebox{0.5\diagramlength}{.}
    \end{tabular}
    \end{equation*}

    The 2-universal property of $A_f$ then induces the multiplication on $A_f$ together with 2-isomorphisms $\psi_x, \psi_y$ filling the following diagram:
    \[\begin{tikzcd}
        {}
        & {A_x \Box A_x}
            \arrow[rr,"m_{x,x,x}"]
            \arrow[rd,"\psi_x",Rightarrow,shorten <=10pt,shorten >=10pt]
        & {}
        & {A_x}
            \arrow[rd,"f_*"]
            \arrow[dd,"\varpi",Rightarrow,shorten <=20pt,shorten >=20pt]
        & {}
        \\{A_f \Box A_f}
            \arrow[ur,"p_x \Box p_x"]
            \arrow[dr,"p_y \Box p_y"']
            \arrow[rr,"m_f",dashed]
        & {}
        & {A_f}
            \arrow[ur,"p_x"']
            \arrow[dr,"p_y"]
        & {}
        & {\mathfrak{M}[x,y]}
        \\ {}
        & {A_y \Box A_y}
            \arrow[rr,"m_{y,y,y}"']
            \arrow[ru,"\psi_y"',Rightarrow,shorten <=10pt,shorten >=10pt]
        & {}
        & {A_y}
            \arrow[ru,"f^*"' ]
        & {}
    \end{tikzcd}.\]

    Similarly, we can provide the unit on $A_f$, which is a 1-morphism $i_f \colon \mathbf{I} \to A_f$ in $\mathfrak{C}$. It suffices to provide the following data:
    \begin{itemize}
        \item 1-morphisms \[\mathbf{I} \xrightarrow{j_x} A_x, \qquad \mathbf{I} \xrightarrow{j_y} A_y;\]
        
        \item one 2-isomorphism given by the following composite:
    \end{itemize}

    \settoheight{\diagramlength}{\includegraphics[height=24mm]{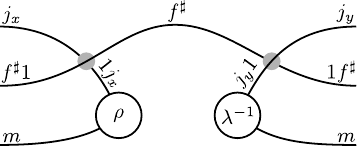}}

    \begin{equation*} \label{eqn:1MorphismUnit}
    \begin{tabular}{@{}ccc@{}}
        \includegraphics[height=24mm]{Pictures/1Morphism/unit.pdf} \raisebox{0.5\diagramlength}{.}
    \end{tabular}
    \end{equation*} 

    The 2-universal property of $A_f$ then induces the unit on $A_f$ together with 2-isomorphisms $\eta_x, \eta_y$ filling the following diagram:
    \[\begin{tikzcd}[sep=30pt]
        {}
        & {A_x}
            \arrow[rd,"f_*"]
            \arrow[dd,"\varpi",Rightarrow,shorten <=20pt,shorten >=20pt,shift left=20pt]
        & {}
        \\ {\mathbf{I}}
            \arrow[ur,"j_x" name=U]
            \arrow[dr,"j_y"' name=D]
            \arrow[r,"i_f",dashed]
        & {A_f}
            \arrow[u,"p_x"']
            \arrow[d,"p_y"]
            \arrow[from=U,Rightarrow,shorten <=5pt,shorten >=5pt,"\eta_x"]
            \arrow[from=D,Rightarrow,shorten <=5pt,shorten >=5pt,"\eta_y"]
        & {\mathfrak{M}[x,y]}
        \\ {}
        & {A_y}
            \arrow[ru,"f^*"']
        & {}
    \end{tikzcd}.\]

    Next, we would like to construct the coherence data on $A_f$: associator $\alpha_f$, left unitor $\lambda_f$ and right unitor $\rho_f$. 

    For the associator, we need to provide a 2-isomorphism \[ \alpha_f \colon m_f \circ (m_f \, \Box \, \mathbf{1}_{A_f}) \to m_f \circ (\mathbf{1}_{A_f} \, \Box \, m_f). \] By the 2-universal property of $A_f$, it suffices to provide two 2-isomorphisms:
    
    \settoheight{\diagramlength}{\includegraphics[height=24mm]{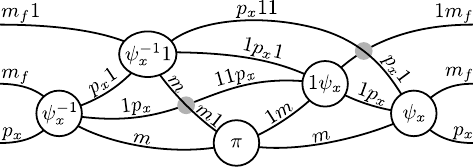}}

    \begin{equation*} \label{eqn:1MorphismAssociator1}
    \begin{tabular}{@{}ccc@{}}
        \includegraphics[height=24mm]{Pictures/1Morphism/associator1.pdf} \raisebox{0.5\diagramlength}{,}
    \end{tabular}
    \end{equation*} 

    \settoheight{\diagramlength}{\includegraphics[height=24mm]{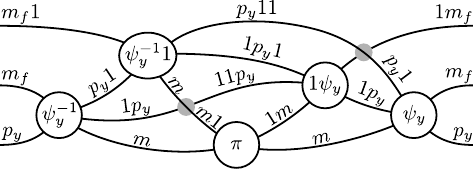}}

    \begin{equation*} \label{eqn:1MorphismAssociator2}
    \begin{tabular}{@{}ccc@{}}
        \includegraphics[height=24mm]{Pictures/1Morphism/associator2.pdf} \raisebox{0.5\diagramlength}{,}
    \end{tabular}
    \end{equation*}

    \noindent and show that the following equation holds:

    \settoheight{\diagramlength}{\includegraphics[height=40mm]{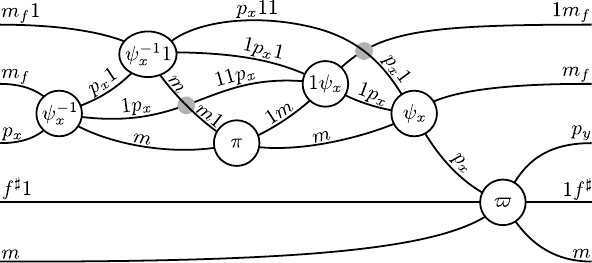}}

    \begin{equation} \label{eqn:1MorphismAssociator3}
    \begin{tabular}{@{}ccc@{}}
        \includegraphics[height=40mm]{Pictures/1Morphism/associator3left.pdf} \raisebox{0.5\diagramlength}{\phantom{.}} \\  \raisebox{5mm}{\phantom{.}} \\ \raisebox{5mm}{\rotatebox{90}{$=$}} \\
        \includegraphics[height=40mm]{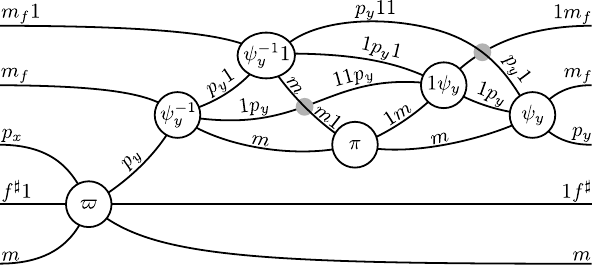} \raisebox{0.5\diagramlength}{.}
    \end{tabular}
    \end{equation}

    For the left unitor, we need to provide a 2-isomorphism \[ \lambda_f \colon m_f \circ (i_f \, \Box \, \mathbf{1}_{A_f}) \to \mathbf{1}_{A_f}. \] By the 2-universal property of $A_f$, it suffices to provide two 2-isomorphisms:

    \settoheight{\diagramlength}{\includegraphics[height=24mm]{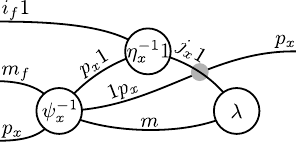}}

    \begin{equation*} \label{eqn:1MorphismLeftUnitor1}
    \begin{tabular}{@{}ccc@{}}
        \includegraphics[height=24mm]{Pictures/1Morphism/leftunitor1.pdf} \raisebox{0.5\diagramlength}{,}
    \end{tabular}
    \end{equation*}

    \settoheight{\diagramlength}{\includegraphics[height=24mm]{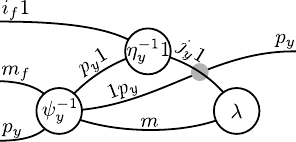}}

    \begin{equation*} \label{eqn:1MorphismLeftUnitor2}
    \begin{tabular}{@{}ccc@{}}
        \includegraphics[height=24mm]{Pictures/1Morphism/leftunitor2.pdf} \raisebox{0.5\diagramlength}{,}
    \end{tabular}
    \end{equation*}
    \noindent and show that the following equation holds:

    \settoheight{\diagramlength}{\includegraphics[height=40mm]{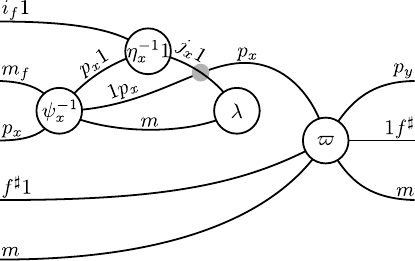}}

    \begin{equation} \label{eqn:1MorphismLeftUnitor3}
    \begin{tabular}{@{}ccc@{}}
        \includegraphics[height=40mm]{Pictures/1Morphism/leftunitor3left.pdf} \raisebox{0.5\diagramlength}{\phantom{.}} \\ \rotatebox{90}{$=$} \\
        \includegraphics[height=40mm]{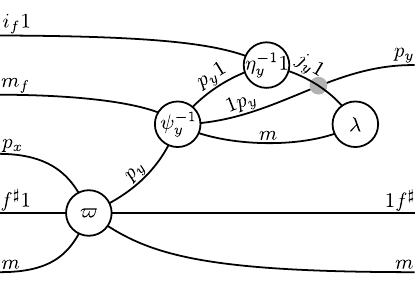} \raisebox{0.5\diagramlength}{.}
    \end{tabular}
    \end{equation}

    For the right unitor, we need to provide a 2-isomorphism \[ \rho_f \colon m_f \circ (\mathbf{1}_{A_f} \, \Box \, i_f) \to \mathbf{1}_{A_f}. \] By the 2-universal property of $A_f$, it suffices to provide two 2-isomorphisms:

    \settoheight{\diagramlength}{\includegraphics[height=24mm]{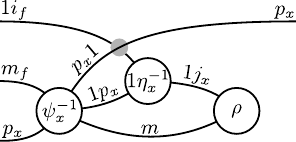}}

    \begin{equation*} \label{eqn:1MorphismRightUnitor1}
    \begin{tabular}{@{}ccc@{}}
        \includegraphics[height=24mm]{Pictures/1Morphism/rightunitor1.pdf} \raisebox{0.5\diagramlength}{,}
    \end{tabular}
    \end{equation*}

    \settoheight{\diagramlength}{\includegraphics[height=24mm]{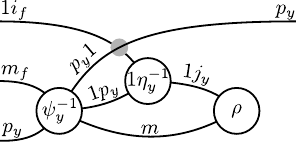}}

    \begin{equation*} \label{eqn:1MorphismRightUnitor2}
    \begin{tabular}{@{}ccc@{}}
        \includegraphics[height=24mm]{Pictures/1Morphism/rightunitor2.pdf} \raisebox{0.5\diagramlength}{,}
    \end{tabular}
    \end{equation*}
    \noindent and show that the following equation holds:

    \settoheight{\diagramlength}{\includegraphics[height=40mm]{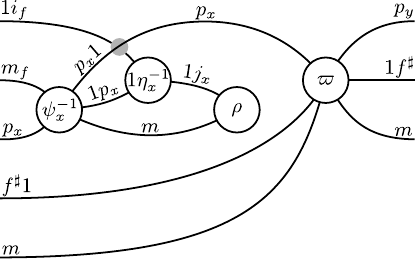}}

    \begin{equation} \label{eqn:1MorphismRightUnitor3}
    \begin{tabular}{@{}ccc@{}}
        \includegraphics[height=40mm]{Pictures/1Morphism/rightunitor3left.pdf} \raisebox{0.5\diagramlength}{\phantom{.}} \\ \rotatebox{90}{$=$} \\
        \includegraphics[height=40mm]{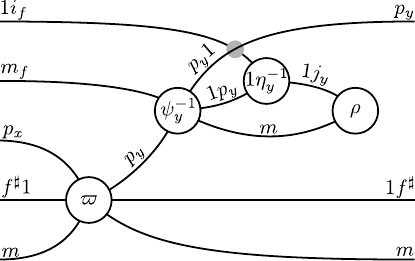} \raisebox{0.5\diagramlength}{.}
    \end{tabular}
    \end{equation}

    Finally, we can check that $(p_x,\psi_x,\eta_x)$ and $(p_y,\psi_y,\eta_y)$ satisfy the coherence equations \eqref{eqn:Algebra1MorphismCoh1}, \eqref{eqn:Algebra1MorphismCoh2}, \eqref{eqn:Algebra1MorphismCoh3}. So they are indeed algebra 1-morphisms. This completes the proof.
\end{proof}

\subsection{On the level of 2-Morphisms} 

Given two 1-morphisms $f,g:x \to y$ in $\mathfrak{M}$ and a 2-morphism $a:f \to g$, we need to assign an algebra $A_a$ in $\mathfrak{C}$ together with algebra 1-morphisms and 2-isomorphisms filling the 2-span.

\[\begin{tikzcd}
    {}
    & {}
    & {A_f}
        \arrow[lldd,"p^f_x"']
        \arrow[rrdd,"p^f_y"]
    & {}
    & {}
    \\ {}
    & {}
        \arrow[dd,Rightarrow,shorten <=15pt,shorten >=15pt,"\phi_x"]
    & {}
    & {}
        \arrow[dd,Rightarrow,shorten <=15pt,shorten >=15pt,"\phi_y"']
    & {}
    \\ {A_x}
    & {}
    & {A_a}
        \arrow[uu,"p_f" description]
        \arrow[dd,"p_g" description]
    & {}
    & {A_y}
    \\ {}
    & {}
    & {}
    & {}
    & {}
    \\ {}
    & {}
    & {A_g}
        \arrow[lluu,"p^g_x"]
        \arrow[rruu,"p^g_y"']
    & {}
    & {}
\end{tikzcd}.\]

\begin{Construction}
    We construct the underlying object of $A_a$ as the following comma object (see Construction \ref{cstr:CommaObject}):
    \[\begin{tikzcd}
        {A_a}
            \arrow[d,"p^a_x"']
            \arrow[r,"p^a_y"]
        & {A_y}
            \arrow[d,"g^*"]
        \\ {A_x} 
            \arrow[r,"f_*"']
            \arrow[ur,Rightarrow,shorten <=10pt,shorten >=10pt,"\varpi_a"]
        & {\mathfrak{M}[x,y]}
    \end{tikzcd}.\] Notably, the 2-morphism $\varpi_a$ is not invertible in general.

    Then, consider the push-forward and pullback along the 2-morphism $a$, we get two composite 2-morphisms of the same shape:
    \[\begin{tikzcd}[sep=30pt]
        {A_f}
            \arrow[d,"p^f_x"']
            \arrow[r,"p^f_y"]
        & {A_y}
            \arrow[d,"f^*"' {name=L}, bend right=30pt]
            \arrow[d,"g^*" {name=R}, bend left=30pt]
            \arrow[from=L,to=R,Rightarrow,"a^*",shorten <=5pt,shorten >=5pt]
        \\ {A_x} 
            \arrow[r,"f_*"']
            \arrow[ur,Rightarrow,shorten <=15pt,shorten >=15pt,"\varpi_f"]
        & {\mathfrak{M}[x,y]}
    \end{tikzcd}, \qquad \begin{tikzcd}[sep=30pt]
        {A_g}
            \arrow[d,"p^g_x"']
            \arrow[r,"p^g_y"]
        & {A_y}
            \arrow[d,"g^*"]
        \\ {A_x} 
            \arrow[r,"g_*" {name=U}, bend left=30pt]
            \arrow[r,"f_*"' {name=D}, bend right=30pt]
            \arrow[ur,Rightarrow,xshift=5pt,yshift=5pt,shorten <=15pt,shorten >=15pt,"\varpi_g"]
            \arrow[from=D,to=U,Rightarrow,xshift=5pt,"a_*"',shorten <=5pt,shorten >=5pt]
        & {\mathfrak{M}[x,y]}
    \end{tikzcd}.\]

    By the 2-universal property of comma object, there are induced 1-morphisms $p_f \colon A_a \to A_f$ and $p_g \colon A_a \to A_g$ together with 2-isomorphisms filling the claimed 2-span above.
\end{Construction}

\begin{Theorem} \label{thm:2SpanOfAlgebrasAssociatedTo2Morphism}
    There is a canonical algebra structure on $A_a$ such that the above 2-span can be promoted to a 2-span of algebras in $\mathfrak{C}$.
\end{Theorem}

\begin{proof}
    The proof is basically the same as the proof of Theorem \ref{thm:AlgebraStructureOn2FiberProduct}, except that we use the 2-universal property of comma object instead of 2-fiber product. The only difference is that we cannot invert the 2-morphism $\varpi_a$ in the construction of multiplication and unit on $A_a$, but this does not affect the rest of the proof.
\end{proof}

\begin{Remark}
    Since the definition of 2-fiber product is symmetric in nature, one can replace the 2-isomorphisms $\varpi_f$ and $\varpi_g$ in the 2-universal property of 2-fiber product by their inverses. However, the definition of comma object is not symmetric so we need to replace the definition of $A_a$ with the following comma object:
    \[\begin{tikzcd}
        {\widetilde{A}_a}
            \arrow[d,"\widetilde{p}^a_x"']
            \arrow[r,"\widetilde{p}^a_y"]
        & {A_y}
            \arrow[d,"f^*"]
            \arrow[dl,Rightarrow,shorten <=10pt,shorten >=10pt,"\widetilde{\varpi}_a"']
        \\ {A_x} 
            \arrow[r,"g_*"']
        & {\mathfrak{M}[x,y]}
    \end{tikzcd},\] in order to make the push-forward and pullback along $a$ compatible.
\end{Remark}

\subsection{Lax Functoriality} 

Suppose $\mathfrak{C}$ has sufficient 2-limits, one can verify that the above assignments on objects, 1-morphisms and 2-morphisms indeed define a lax 3-functor from $\mathfrak{M}$ to $\mathbf{Span}_{(3,2)}(\mathbf{Alg}(\mathfrak{C}))$.

First, notice that for each object $x$ in $\mathfrak{M}$, the identity 1-morphism $\mathbf{1}_x \colon x \to x$ is assigned to the identity span on $A_x$ canonically. This strictifies the higher unitors for the lax 3-functor structure, i.e., this is a \emph{normalized} lax 3-functor.

\begin{Theorem} \label{thm:Laxator1Morphism}
    Given a pair of composable 1-morphisms $x \xrightarrow{f} y \xrightarrow{g} z$ in $\mathfrak{M}$, there is a canonical 1-morphism \[\Phi_{f,g} \colon A_{g \circ f} \to A_f \times_{A_y} A_g.\]
\end{Theorem}

\begin{proof}
    Consider the following diagram:
    \[\begin{tikzcd}
        {A_f \times_{A_y} A_g}
            \arrow[d]
            \arrow[r]
        & {A_g}
            \arrow[r,"p^g_z"]
            \arrow[d,"p^g_y"']
        & {A_z}
            \arrow[d,"g^*"]
        \\ {A_f}
            \arrow[r,"p^f_y"]
            \arrow[d,"p^f_x"']
            \arrow[ur,Rightarrow,shorten <=20pt,shorten >=20pt]
        & {A_y}
            \arrow[r,"g^*"]
            \arrow[d,"f^*"']
            \arrow[ur,Rightarrow,shorten <=20pt,shorten >=20pt]
        & {\mathfrak{M}[y,z]}
            \arrow[d,"f^*"]
        \\ {A_x}
            \arrow[r,"f_*"']
            \arrow[ur,Rightarrow,shorten <=20pt,shorten >=20pt]
        & {\mathfrak{M}[x,y]}
            \arrow[r,"g_*"']
            \arrow[ur,Rightarrow,shorten <=10pt,shorten >=10pt]
        & {\mathfrak{M}[x,z]}
    \end{tikzcd},\] where the bottom right square is filled by the interchanger between $f$ and $g$ in $\mathfrak{M}$, and the other three squares are filled by the canonical 2-isomorphisms for the 2-fiber products. Notice that the composition of the bottom row is the push-forward $(g \circ f)_*$ and the composition of the right column is the pullback $(g \circ f)^*$. Hence, $\Phi_{f,g}$ is induced by the 2-universal property of the 2-fiber product $A_{g \circ f}$.
\end{proof}

\begin{Remark}
    The 1-morphism $\Phi_{f,g}$ constructed above is not an isomorphism in general. The failure is caused by the fact that the bottom right square in the above diagram is not a 2-fiber product square in general.
\end{Remark}

\begin{Theorem} \label{thm:Laxator2Morphism}
    Given composable 1-morphisms $x \xrightarrow{f} y \xrightarrow{g} z \xrightarrow{h} w$ in $\mathfrak{M}$, there is a canonical 2-isomorphism filling the following diagram:
    \[\begin{tikzcd}
        {A_{h \circ g \circ f}}
            \arrow[r,"\Phi_{f,hg}"]
            \arrow[d,"\Phi_{gf,h}"']
        & {A_{f} \times_{A_y} A_{h \circ g}}
            \arrow[d,"1 \Phi_{g,h}"]
        \\ {A_{g \circ f} \times_{A_z} A_{h}}
            \arrow[r,"\Phi_{f,g} 1"']
            \arrow[ur,Rightarrow,shorten <=10pt,shorten >=10pt,"\Phi_{f,g,h}"]
        &{A_f \times_{A_y} A_{g} \times_{A_z} A_{h}}
    \end{tikzcd}.\]
\end{Theorem}

\begin{proof}
    Consider the following diagram:
    \[\begin{tikzcd}
        {A_f \times_{A_y} A_g \times_{A_z} A_h}
            \arrow[r]
            \arrow[d]
        & {A_g \times_{A_z} A_h}
            \arrow[r,color=blue]
            \arrow[d,color=blue]
        & {A_h}
            \arrow[r,color=blue]
            \arrow[d]
        & {A_w}
            \arrow[d,"h^*",color=blue]
        \\ {A_f \times_{A_y} A_g}
            \arrow[r,color=red]
            \arrow[d,color=red]
        & {A_g}
            \arrow[r,color=red]
            \arrow[d,color=blue]
        & {A_z}
            \arrow[r,"h_*"']
            \arrow[d,"g^*",color=red]
        & {\mathfrak{M}[z,w]}
            \arrow[d,"g^*",color=blue]
        \\ {A_f}
            \arrow[r]
            \arrow[d,color=red]
        & {A_y}
            \arrow[r,"g_*"',color=blue]
            \arrow[d,"f^*"]
        & {\mathfrak{M}[y,z]}
            \arrow[r,"h_*"',color=blue]
            \arrow[d,"f^*",color=red]
        & {\mathfrak{M}[y,w]}
            \arrow[d,"f^*"]
        \\ {A_x}
            \arrow[r,"f_*"',color=red]
        & {\mathfrak{M}[x,y]}
            \arrow[r,"g_*"',color=red]
        & {\mathfrak{M}[x,z]}
            \arrow[r,"h_*"']
        & {\mathfrak{M}[x,w]}
    \end{tikzcd},\] where three bottom right squares are filled by the interchangers, and the other squares are filled by the canonical 2-isomorphisms for the 2-fiber products. In comparison with the proof of Theorem \ref{thm:Laxator1Morphism}, one can first consider the bottom left four squares painted in red, where by the 2-universal property of $A_{f} \times_{A_y} A_g$, the above diagram factors via $\Phi_{f,g}$ through the following four squares:
    \[\begin{tikzcd}
        {A_{g \circ f,h} \times_{A_z} A_h}
            \arrow[r]
            \arrow[d]
        & {A_{h}}
            \arrow[r]
            \arrow[d]
        & {A_w}
            \arrow[d,"h^*"]
        \\ {A_{g \circ f}}
            \arrow[r]
            \arrow[d]
        & {A_z}
            \arrow[r,"h_*"']
            \arrow[d,"(g \circ f)^*"]
        & {\mathfrak{M}[z,w]}
            \arrow[d,"(g \circ f)^*"]
        \\ {A_x}
            \arrow[r,"(g \circ f)_*"']
        & {\mathfrak{M}[x,z]}
            \arrow[r,"h_*"']
        & {\mathfrak{M}[x,w]}
    \end{tikzcd}.\] Again, the 2-universal property of $A_f \times_{A_y} A_g \times_{A_z} A_h$ further factors the outer square via $\Phi_{gf,h}$ through the square:
    \[\begin{tikzcd}
        {A_{h \circ g \circ f}}
            \arrow[r]
            \arrow[d]
        & {A_w}
            \arrow[d,"(h \circ g \circ f)^*"]
        \\ {A_{x}}
            \arrow[r,"(h \circ g \circ f)_*"']
        & {\mathfrak{M}[x,w]}
    \end{tikzcd}.\] This completes the construction of the source of the desired 2-isomorphism $\Phi_{f,g,h}$. The target can be constructed similarly by first considering the top right four squares painted in blue. Finally, one can check that the two ways of factoring the outer square agree via a canonical 2-isomorphism by the 2-universal property of the triple 2-fiber product.
\end{proof}

\begin{Corollary} 
    Generalizing Theorem \ref{thm:Laxator1Morphism} and Theorem \ref{thm:Laxator2Morphism} to a sequence of four composable 1-morphisms in $\mathfrak{M}$, one can verify that the above assignments indeed define a normalized lax 3-functor from $\mathfrak{M}$ to $\mathbf{Span}_{(3,2)}(\mathbf{Alg}(\mathfrak{C}))$.
\end{Corollary}

Lastly, recall from Definition \ref{def:ModulesInEnriched2Category} that for any algebra $A$ in $\mathfrak{C}$ and any $\mathfrak{C}$-enriched 2-category $\mathfrak{M}$, the 2-category of modules $\mathbf{Mod}_{\mathfrak{M}}(A)$ is defined to be the 2-category of $\mathfrak{C}$-enriched 2-functors: $\mathbf{Cat}(\mathfrak{C})(\mathbb{B} A, \mathfrak{M})$. 

\begin{Theorem} \label{thm:Embed2CatOfModulesInto2SpansOfAlgebras}
    Combining the above lax 3-functor and Proposition \ref{prop:DeloopingOfAlgebrasIs3Functor}, for $\mathfrak{C}$ a monoidal 2-category with 2-fiber products and comma objects, we obtain a lax 3-functor, which takes an $A$-module in a $\mathfrak{C}$-enriched 2-category $\mathfrak{M}$ to a 2-span of algebras in $\mathfrak{C}$ under $A$:
\[\mathbf{Mod}_{\mathfrak{M}}(A) \to \mathbf{Span}_{(3,2)}({}^{A/}\mathbf{Alg}(\mathfrak{C})).\]
\end{Theorem}

\begin{Remark}
    The significance of Theorem \ref{thm:Embed2CatOfModulesInto2SpansOfAlgebras} lies in its ability to provide a high-level justification for practical intuitions. While the theorem transforms the familiar notion of modules into the less familiar setting of 2-spans of algebras via a lax 3-functor that does not fully preserve structures, this shift is motivated by practical considerations. In computations, modules often need to be analyzed in terms of algebras, and this framework formalizes and supports that reduction, bridging theoretical abstraction with computational utility.
\end{Remark}

%% file: Applications.tex
\section{Applications and Further Discussions} 

\subsection{Applications to Module Categories} \label{sec:ApplicationToCat}

Once we have established the general theory of modules as spans of algebras, we can apply it to various examples of monoidal 2-categories.

    The first motivating example, as previously discussed in the title and the introduction, is when $\mathfrak{C} = \mathbf{Cat}$ is the 2-category of (small) categories, with the monoidal structure given by the Cartesian product. In this case, algebras in $\mathfrak{C}$ are precisely monoidal categories, and modules in $\mathfrak{C}$ are module categories over monoidal categories.

    For this particular example, let us spell out the details of the embedding lax 3-functor from Theorem \ref{thm:Embed2CatOfModulesInto2SpansOfAlgebras}. Given a monoidal category $\mathcal{C}$, a left $\mathcal{C}$-module category consists of a category $\mathcal{M}$ together with a monoidal functor \[F^\mathcal{M}: \mathcal{C} \to \mathbf{End}(\mathcal{M}),\] where $\mathbf{End}(\mathcal{M})$ is the monoidal category of endofunctors on $\mathcal{M}$ with composition as the monoidal structure.

    Given two left $\mathcal{C}$-module categories $\mathcal{M}$ and $\mathcal{N}$, a functor $F:\mathcal{M} \to \mathcal{N}$ admits a left $\mathcal{C}$-module functor structure if and only if there exists a monoidal functor \[\Phi^F: \mathcal{C} \to \mathbf{End}(\mathcal{M}) \times^{F} \mathbf{End}(\mathcal{N}),\] where $\mathbf{End}(\mathcal{M}) \times^{F} \mathbf{End}(\mathcal{N})$ is the 2-fiber product of categories: 
    \[\begin{tikzcd}[sep=large]
    {\mathbf{End}(\mathcal{M}) \times^F \mathbf{End}(\mathcal{N})}
        \arrow[d,"\mathbf{pr2}"']
        \arrow[r,"\mathbf{pr1}"]
    & {\mathbf{End}(\mathcal{M})}
        \arrow[d,"F_*"]
        \arrow[dl,Rightarrow,shorten <=25pt, shorten >=20pt,"\phi^F"']
    \\ {\mathbf{End}(\mathcal{N})}
        \arrow[r,"F^*"']
    & {\mathbf{Fun}(\mathcal{M},\mathcal{N})}
    \end{tikzcd},\] where $F_*$ is the post-composition with $F$ and $F^*$ is the pre-composition with $F$. 

    An object in $\mathbf{End}(\mathcal{M}) \times^F \mathbf{End}(\mathcal{N})$ is a triple $(P,Q,\xi)$, where:
    \begin{itemize}
        \item $P$ is an endofunctor on $\mathcal{M}$,
        
        \item $Q$ is an endofunctor on $\mathcal{N}$,
        
        \item $\xi:F \circ_1 P \to Q \circ_1 F$ is a natural isomorphism.
    \end{itemize}
    
    A morphism in $\mathbf{End}(\mathcal{M}) \times^F \mathbf{End}(\mathcal{N})$ from $(P_0,Q_0,\xi_0)$ to $(P_1,Q_1,\xi_1)$ consists of a pair of natural transformations $\mu:P_0 \to P_1$ and $\nu:Q_0 \to Q_1$ such that $\xi_1 \circ_2 (F \circ_1 \mu) = (\nu \circ_1 F) \circ_2 \xi_0$. Composition of morphisms in $\mathbf{End}(\mathcal{M}) \times^F \mathbf{End}(\mathcal{N})$ is given by the vertical composition of natural transformations.
    
    Lastly, given two left $\mathcal{C}$-module functors $F,G:\mathcal{M} \to \mathcal{N}$, a natural transformation $\varphi:F \to G$ is a left $\mathcal{C}$-module natural transformation if and only if the monoidal functors \[\Phi^F: \mathcal{C} \to \mathbf{End}(\mathcal{M}) \times^{F} \mathbf{End}(\mathcal{N})\] and \[\Phi^G: \mathcal{C} \to \mathbf{End}(\mathcal{M}) \times^{G} \mathbf{End}(\mathcal{N})\] factor through a monoidal functor \[\mathcal{C} \to \mathbf{End}(\mathcal{M}) \times^{\varphi} \mathbf{End}(\mathcal{N}).\]

    Here, $\mathbf{End}(\mathcal{M}) \times^{\varphi} \mathbf{End}(\mathcal{N})$ is defined as the 2-fiber product of categories:
    \[\begin{tikzcd}
        {\mathbf{End}(\mathcal{M}) \times^{\varphi} \mathbf{End}(\mathcal{N})}
            \arrow[d,"\mathbf{pr2}"']
            \arrow[r,"\mathbf{pr1}"]
        & {\mathbf{End}(\mathcal{M}) \times^F \mathbf{End}(\mathcal{N})}
            \arrow[d,"(\varphi^*)_*"]
            \arrow[dl,Rightarrow,shorten <=25pt, shorten >=20pt,"\phi^\varphi"']
        \\ {\mathbf{End}(\mathcal{M}) \times^G \mathbf{End}(\mathcal{N})}
            \arrow[r,"(\varphi_*)^*"']
        & {\mathbf{End}(\mathcal{M}) {}_{G^*} ;_{F_*} \mathbf{End}(\mathcal{N})}
    \end{tikzcd},\] where:
    \begin{itemize}
        \item The natural transformation $\varphi^* \colon F^* \to G^*$ is induced by the functoriality of pullback along functors, which further induces the natural transformation
        \[\mathbf{End}(\mathcal{M}) \times^{F} \mathbf{End}(\mathcal{N}) \xrightarrow{(\varphi^*)_*} \mathbf{End}(\mathcal{M}) \times^{G^*}_{F_*} \mathbf{End}(\mathcal{N}),\]
        \[(P,Q,F \circ_1 P \xrightarrow{\xi} Q \circ_1 F) \longmapsto (P,Q,F \circ_1 P \xrightarrow{\xi} Q \circ_1 F \xrightarrow{Q \circ_1 \varphi} Q \circ_1 G);\]

        \item The natural transformation $\varphi_*:F_* \to G_*$ is induced by the functoriality of push-forward along functors, which further induces the natural transformation
        \[\mathbf{End}(\mathcal{M}) \times^{G} \mathbf{End}(\mathcal{N}) \xrightarrow{(\varphi_*)^*} \mathbf{End}(\mathcal{M}) \times^{G^*}_{F_*} \mathbf{End}(\mathcal{N}),\]
        \[(P,Q,G \circ_1 P \xrightarrow{\xi} Q \circ_1 G) \longmapsto (P,Q,F \circ_1 P \xrightarrow{\varphi \circ_1 P} G \circ_1 P \xrightarrow{\xi} Q \circ_1 G);\]

        \item $\mathbf{End}(\mathcal{M}) {}_{G^*} ;_{F_*} \mathbf{End}(\mathcal{N})$ is the comma category:
        \[\begin{tikzcd}
        {\mathbf{End}(\mathcal{M}) \times^{G^*}_{F_*} \mathbf{End}(\mathcal{N})}
            \arrow[d,"\mathbf{pr2}"']
            \arrow[r,"\mathbf{pr1}"]
        & {\mathbf{End}(\mathcal{M})}
            \arrow[d,"F_*"]
            \arrow[dl,Rightarrow,shorten <=25pt, shorten >=20pt]
        \\ {\mathbf{End}(\mathcal{N})}
            \arrow[r,"G^*"']
        & {\mathbf{Fun}(\mathcal{M},\mathcal{N})}
        \end{tikzcd}.\]
    \end{itemize}

    An object in $\mathbf{End}(\mathcal{M}) \times^\varphi \mathbf{End}(\mathcal{N})$ consists of $(P,Q,\xi^F,\xi^G)$ where
    \begin{itemize}
        \item $P$ is an endofunctor on $\mathcal{M}$,

        \item $Q$ is an endofunctor on $\mathcal{N}$,

        \item $\xi^F: F \circ_1 P \to Q \circ_1 F$ is a natural isomorphism,

        \item $\xi^G: G \circ_1 P \to Q \circ_1 G$ is a natural isomorphism,
    \end{itemize} such that $\xi^G \circ_2 (\varphi \circ_1 P) = (Q \circ_1 \varphi) \circ_2 \xi^F$.
    
    A morphism from $(P_0,Q_0,\xi^F_0,\xi^G_0)$ to $(P_1,Q_1,\xi^F_1,\xi^G_1)$ consists of a pair of natural transformations $\mu:P_0 \to P_1$ and $\nu:Q_0 \to Q_1$ such that $\xi^F_1 \circ_2 (F \circ_1 \mu) = (\nu \circ_1 F) \circ_2 \xi^F_0$ and $\xi^G_1 \circ_2 (G \circ_1 \mu) = (\nu \circ_1 G) \circ_2 \xi^G_0$.
    
    Composition of morphisms in $\mathbf{End}(\mathcal{M}) \times^\varphi \mathbf{End}(\mathcal{N})$ is given by the vertical composition of natural transformations.

    In summary, a left $\mathcal{C}$-module natural transformation \[\begin{tikzcd}[column sep=40pt]
        {\mathcal{M}}
            \arrow[r,bend left=50pt,"F"{name=U}]
            \arrow[r,bend right=50pt,"G"'{name=D}]
        & {\mathcal{N}}
        \arrow[Rightarrow, from=U, to=D, "\varphi", shorten <=10pt, shorten >=10pt]
    \end{tikzcd} \] is sent to the following 2-span of monoidal categories under $\mathcal{C}$:
    \[ \begin{tikzcd}[column sep=40pt]
        {}
        & {\mathbf{End}(\mathcal{M}) \times^{F} \mathbf{End}(\mathcal{N})}
            \arrow[dl]
            \arrow[dr]
        & {}
        \\ {\mathbf{End}(\mathcal{M})}
        & {\mathbf{End}(\mathcal{M}) \times^{\varphi} \mathbf{End}(\mathcal{N})}
            \arrow[d]
            \arrow[u]
        & {\mathbf{End}(\mathcal{N})}
        \\ {}
        & {\mathbf{End}(\mathcal{M}) \times^{G} \mathbf{End}(\mathcal{N})}
            \arrow[ur]
            \arrow[ul]
        & {}
    \end{tikzcd}.\]

\subsection{Examples of Algebras and Modules in Monoidal 2-Categories} \label{sec:OtherExamples}

In the same spirit, we view modules in other monoidal 2-categories as spans of algebras. Here we list some notable examples of monoidal 2-categories satisfying the assumptions of Theorem \ref{thm:Embed2CatOfModulesInto2SpansOfAlgebras}, along with their corresponding notions of algebras and modules.

\begin{Example}
    Suppose $\mathfrak{C} = \mathbf{2Vect}$ is the 2-category of linear categories, linear functors and natural transformations, with the monoidal structure given by the Deligne-Kelly tensor product \cite{Del07}. Algebras in $\mathfrak{C}$ are linear monoidal categories, and modules in $\mathfrak{C}$ are module categories over linear monoidal categories \cite{Ost1}.
\end{Example}

\begin{Remark}
    The 2-category of all small linear categories clearly admits 2-fiber products and comma objects. However, while the standard Deligne-Kelly tensor product exists in this setting, it is ill-suited for representation theory. Specifically, the tensor product functor is not right exact in each variable, which breaks the compatibility with Abelian structures. To properly incorporate homological algebra, one can adopt the more sophisticated framework of \emph{presentable stable linear categories} \cite{Lur17}.
\end{Remark}

\begin{Remark}
    We caution the reader that it is insufficient to work exclusively within the 2-category \emph{finite semisimple linear categories}. On one hand, the objects of interest—--finite semisimple module categories over fusion categories—--are well-understood in terms of internal algebras \cite{Ost1} and Morita theory \cite{ENO02}. On the other hand, our enriched perspective constructs 2-spans using 2-fiber products and comma objects, and $\mathbf{2Vect}$ is not closed under these 2-limits.

    For a counterexample regarding finiteness, consider the 2-fiber product of the diagonal functor $\Delta \colon \mathbf{Vect} \to \mathbf{Vect} \times \mathbf{Vect}$ with itself. An object in this category corresponds to a finite dimensional vector space $V$ equipped with an automorphism (essentially two isomorphisms $T_0, T_1 \colon V \to W$ which reduce to $T_1^{-1} \circ T_0$). This is equivalent to the category of finite dimensional representations of $\mathbb{Z}$, which is not a finite category.

    For a counterexample regarding semisimplicity, consider the comma category of the identity functor $\text{Id} \colon \mathbf{Vect} \to \mathbf{Vect}$ with itself. An object here consists of two vector spaces $V, W$ and a linear map $f: V \to W$. This is equivalent to the category of finite dimensional representations of the quiver $A_2$ ($\bullet \to \bullet$), which admits non-trivial extensions and is therefore not semisimple.
\end{Remark}

\begin{Example}
    Fix a finite group $G$. If $\mathfrak{C} = \mathbf{2Rep}(G)$ is the 2-category of linear categories with $G$-actions \cite[Construction 2.1.12 for the finite semisimple version]{DR}, with the monoidal structure induced by the Deligne tensor product of the underlying linear categories, then:
    \begin{itemize}
        \item An algebra in $\mathfrak{C}$ consists of a linear monoidal category $\mathcal{C}$ together with a $G$-action on $\mathcal{C}$ via monoidal auto-equivalences, i.e. a 2-group morphism $G \to \mathbf{Aut}_\otimes(\mathcal{C})$;
        
        \item A module in $\mathfrak{C}$ consists of a module category $\mathcal{M}$ with $G$-equivariant structure on the module category data.
    \end{itemize}
\end{Example}

\begin{Example}
    Fix a finite group $G$ and a group 4-cocycle $\pi$. Let $\mathfrak{C} = \mathbf{2Vect}^\pi_G$ be the 2-category of $G$-graded linear categories \cite[Construction 2.1.1 for the finite semisimple version]{DR}, with the monoidal structure induced by convolution on $G$ and pentagonator twisted by $\pi$. 
    
    Then algebras in $\mathfrak{C}$ are $G$-graded linear monoidal categories\footnote{Here $G$-gradings are not assumed to be faithful, i.e. it is possible to have $G$-graded components which are zero.} on which the anomaly $\pi$ is restricted to an exact cocycle, and modules in $\mathfrak{C}$ are $G$-graded module categories, see \cite[Example 1.5.5 for the finite semisimple version]{D7}.
\end{Example}

\begin{Example} \label{exmp:CentralModuleFusionCategories}
    Take a braided fusion category $\mathcal{B}$. Suppose $\mathfrak{C} = \mathbf{Mod}(\mathcal{B})$ is the 2-category of $B$-module categories \cite[Construction 2.1.19 for the finite semisimple version]{DR}, with the monoidal structure induced by the relative Deligne tensor product $\boxtimes_\mathcal{B}$. An algebra in $\mathfrak{C}$ consists of a central $\mathcal{B}$-module monoidal category $\mathcal{C}$, i.e., a monoidal category $\mathcal{C}$ together with a braided functor $\mathcal{B} \to \mathcal{Z}_1(\mathcal{C})$, see \cite[Example 1.5.6 for the finite semisimple version]{D7}. 
    
    A left module over $(\mathcal{C},\mathcal{B} \xrightarrow{F} \mathcal{Z}_1(\mathcal{C}))$ is an ordinary left $\mathcal{C}$-module category $\mathcal{M}$. It is equipped with an induced left $\mathcal{B}$-action on $\mathcal{M}$ by \[\mathcal{B} \xrightarrow{F} \mathcal{Z}_1(\mathcal{C}) \xrightarrow{\text{forget}} \mathcal{C} \to \mathbf{End}(\mathcal{M}).\]
\end{Example}

\subsection{Scenarios Beyond the Enriched 2-Categorical Framework} \label{sec:BeyondEnriched2Categories}
There are other intriguing examples of monoidal 2-categories where the conditions of Theorem \ref{thm:Embed2CatOfModulesInto2SpansOfAlgebras} do not hold. In this section, we offer additional insights on extending the main theorem to encompass these cases.

Consider the case when $\mathfrak{C} = \mathbf{MCat}$ is the 2-category of monoidal categories with the Cartesian product. First, we notice that this monoidal 2-category is not closed, as the internal hom of two monoidal categories is not even defined in general.

\begin{Remark}
    Algebras in $\mathfrak{C}$ are braided monoidal categories, or equivalently a category equipped with two compatible monoidal structures \cite{BFSV03}. This fact generalizes the Eckmann-Hilton arguments for monoids: a set with two compatible monoid structures is exactly a commutative monoid.
\end{Remark}

\medskip \noindent \textit{On the level of Objects.} Given a braided monoidal category $\mathcal{A}$, a left module of $\mathcal{A}$ in $\mathfrak{C}$ consists of an underlying monoidal category $\mathcal{M}$ and a braided functor $F^\mathcal{M}:\mathcal{A} \to \mathcal{Z}_1(\mathcal{M})$, or a \textit{central $\mathcal{A}$-module monoidal category}, as discussed in Example \ref{exmp:CentralModuleFusionCategories}.

\begin{Definition}
    For the precise definition of morphisms between central module monoidal categories, we follow Henriques, Penneys, and Tener \cite{HPT23}. Let us fix a braided monoidal category $\mathcal{A}$. 
    
    Given two central $\mathcal{A}$-module monoidal categories $\mathcal{M}$ and $\mathcal{N}$ in $\mathfrak{C}$, a \emph{central $\mathcal{A}$-module monoidal functor} consists of a monoidal functor $G:\mathcal{M} \to \mathcal{N}$ together with a monoidal natural isomorphism\footnote{These structures encompass the first two coherence conditions in \cite[Definition 3.2]{HPT23}.} \[\begin{tikzcd}[column sep=large]
    {}
    & {\mathcal{Z}_1(\mathcal{M})}
        \arrow[r,"\mathrm{forget}"]
        \arrow[dd,Rightarrow,shorten <=10pt,shorten >=10pt,"\psi^G"]
    & {\mathcal{M}}
        \arrow[dd,"G"]
    \\  {\mathcal{A}}
        \arrow[ur,"F^\mathcal{M}"]
        \arrow[dr,"F^\mathcal{N}"']
    & {}
    & {}
    \\ {}
    & {\mathcal{Z}_1(\mathcal{N})}
        \arrow[r,"\mathrm{forget}"']
    & {\mathcal{N}}
\end{tikzcd},\] subject to an additional coherence condition, which is denoted by (7) in \cite[Definition 3.2]{HPT23}.

Given two central $\mathcal{A}$-module monoidal functors $G,H \colon \mathcal{M} \to \mathcal{N}$, a \emph{central $\mathcal{A}$-module monoidal natural transformation} consists of a monoidal natural transformation $\varphi \colon G \to H$ such that 
\[\begin{tikzcd}
    {}
    & {\mathcal{Z}_1(\mathcal{M})}
        \arrow[r]
        \arrow[dd,Rightarrow,shorten <=10pt,shorten >=10pt,"\psi^G"]
    & {\mathcal{M}}
        \arrow[dd,bend left=35pt,"H"]
        \arrow[dd,bend right=35pt,"G"']
    \\  {\mathcal{A}}
        \arrow[ur,"F^\mathcal{M}"]
        \arrow[dr,"F^\mathcal{N}"']
    & {}
    & {}
        \arrow[l,Leftarrow,xshift=25pt,shorten <=15pt,shorten >=20pt,"\varphi"']
    \\ {}
    & {\mathcal{Z}_1(\mathcal{N})}
        \arrow[r]
    & {\mathcal{N}}
\end{tikzcd} \quad = \quad \begin{tikzcd}
    {}
    & {\mathcal{Z}_1(\mathcal{M})}
        \arrow[r]
        \arrow[dd,Rightarrow,shorten <=10pt,shorten >=10pt,"\psi^H"]
    & {\mathcal{M}}
        \arrow[dd,"H"]
    \\  {\mathcal{A}}
        \arrow[ur,"F^\mathcal{M}"]
        \arrow[dr,"F^\mathcal{N}"']
    & {}
    & {}
    \\ {}
    & {\mathcal{Z}_1(\mathcal{N})}
        \arrow[r]
    & {\mathcal{N}}
\end{tikzcd},\] which corresponds to the condition (9) in \cite[Definition 3.3]{HPT23}.
\end{Definition}

We would like to interpret central $\mathcal{A}$-module monoidal functors and natural transformations into spans of braided monoidal categories over $\mathcal{A}$.

\medskip \noindent \textit{On the level of 1-Morphisms.} Consider the 2-fiber product of monoidal categories: \[\begin{tikzcd}
    {\mathcal{Z}_1(\mathcal{M}) \times_{\mathcal{Z}_1(G)} \mathcal{Z}_1(\mathcal{N})}
        \arrow[r,"\mathbf{pr1}"]
        \arrow[d,"\mathbf{pr2}"']
    & {\mathcal{Z}_1(\mathcal{M})}
        \arrow[d,"G_*"]
        \arrow[dl,Rightarrow,shorten <=25pt, shorten >=20pt]
    \\ {\mathcal{Z}_1(\mathcal{N})}
        \arrow[r,"G^*"']
    & {\mathcal{Z}_1(G)}
\end{tikzcd}\] where $\mathcal{Z}_1(G)$ is the monoidal centralizer, in the sense of \cite[§5.3.1]{Lur17}. In the context of tensor categories, the monoidal centralizer of a subcategory appeared in \cite[§3.6]{DGNO}.

\begin{Remark}
    The monoidal centralizer of a monoidal functor $G \colon \mathcal{M} \to \mathcal{N}$ is equivalent to $\mathbf{End}_{\mathcal{M} \boxtimes \mathcal{N}^{mp}}({}_{\langle F \rangle} \mathcal{N})$ and $\mathbf{Fun}_{\mathcal{M} \boxtimes \mathcal{M}^{mp}}(\mathcal{M},{}_{\langle F \rangle} \mathcal{N}_{\langle F \rangle})$, where the subscripts ${}_{\langle F \rangle}$ indicate that we are considering the module category structure induced by $F$.
\end{Remark}

\begin{Remark} \label{rmk:CoherenceForMonoidalCentralizer}
    A generic object in this 2-fiber product consists of a half-braiding $(m,\beta^m)$ in $\mathcal{M}$, a half-braiding $(n,\beta^n)$ in $\mathcal{N}$ and an invertible 2-morphism $\xi:G(m) \to n$ in $\mathcal{N}$ compatible with half-braidings, i.e., the following diagram commutes: \[\begin{tikzcd}[sep=large]
    {G(m \otimes_\mathcal{M} z)}
        \arrow[r,"G(\beta^m_z)"]
    & {G(z \otimes_\mathcal{M} m)}
    \\ {G(m) \otimes_\mathcal{N} G(z)}
        \arrow[d,"\xi \otimes G(z)"']
        \arrow[u,"\gamma^G_{m,z}"]
    & {G(z) \otimes_\mathcal{N} G(m)}
        \arrow[d,"G(z) \otimes \xi"]
        \arrow[u,"\gamma^G_{z,m}"']
    \\ {n \otimes_\mathcal{N} G(z)}
        \arrow[r,"\beta^n_{Gz}"']
    & {G(z) \otimes_\mathcal{N} n}
\end{tikzcd}\] for any object $z$ in $\mathcal{M}$. This coherence condition is precisely the coherence condition (7) in \cite[Definition 3.2]{HPT23} mentioned earlier.
\end{Remark}

\begin{Construction}
    Monoidal category $\mathcal{Z}_1(\mathcal{M}) \times_{\mathcal{Z}_1(G)} \mathcal{Z}_1(\mathcal{N})$ is equipped with a braiding, given on objects $(m_0,n_0,\xi_0)$ and $(m_1,n_1,\xi_1)$ by \[\left( m_0 \otimes_\mathcal{M} m_1 \xrightarrow{\beta^{m_0}_{m_1}} m_1 \otimes_\mathcal{M} m_0, \quad n_0 \otimes_\mathcal{N} n_1 \xrightarrow{\beta^{n_0}_{n_1}} n_1 \otimes_\mathcal{N} n_0 \right).\]
\end{Construction}

\begin{Proposition}
    For any monoidal functor $G \colon \mathcal{M} \to \mathcal{N}$, a central $\mathcal{A}$-module monoidal functor structure $\psi^G$ corresponds to a braided functor \[\mathcal{A} \to \mathcal{Z}_1(\mathcal{M}) \times_{\mathcal{Z}_1(G)} \mathcal{Z}_1(\mathcal{N}); \quad x \mapsto \left( F^\mathcal{M}(x), F^\mathcal{N}(x), G(F^\mathcal{M}(x)) \xrightarrow{\psi^G_x} F^\mathcal{N}(x) \right).\]
\end{Proposition}

\begin{proof}
    We can rewrite the data of central $\mathcal{A}$-module monoidal functor into a monoidal natural isomorphism \[\begin{tikzcd}[sep=large]
        {\mathcal{A}}
            \arrow[r,"F^\mathcal{M}"]
            \arrow[d,"F^\mathcal{N}"']
        & {\mathcal{Z}_1(\mathcal{M})}
            \arrow[d,"G_*"]
            \arrow[dl,Rightarrow,shorten <=20pt, shorten >=15pt,"\phi^G"']
        \\ {\mathcal{Z}_1(\mathcal{N})}
            \arrow[r,"G^*"']
        & {\mathcal{Z}_1(G)}
    \end{tikzcd},\] satisfying a coherence condition corresponding to (7) in \cite[Definition 3.2]{HPT23}. By the 2-universal property of 2-fiber product, this induces the monoidal functor $\mathcal{A} \to \mathcal{Z}_1(\mathcal{M}) \times_{\mathcal{Z}_1(G)} \mathcal{Z}_1(\mathcal{N})$. The coherence condition ensures that this monoidal functor preserves the braidings.
\end{proof}

\noindent \textit{On the level of 2-morphisms.} Given two central $\mathcal{A}$-module monoidal functors $(G,\psi^G)$ and $(H,\psi^H)$, for a monoidal natural transformation $\varphi:G \to H$, consider the 2-fiber product of monoidal categories: \[\begin{tikzcd}
    {\mathcal{Z}_1(\mathcal{M}) \times^\varphi \mathcal{Z}_1(\mathcal{N})}
        \arrow[d,"\mathbf{pr2}"']
        \arrow[r,"\mathbf{pr1}"]
    & {\mathcal{Z}_1(\mathcal{M}) \times_{\mathcal{Z}_1(G)} \mathcal{Z}_1(\mathcal{N})}
        \arrow[d,"(\varphi_*)^*"]
        \arrow[dl,Rightarrow,shorten <=25pt, shorten >=20pt]
    \\ {\mathcal{Z}_1(\mathcal{M}) \times_{\mathcal{Z}_1(H)} \mathcal{Z}_1(\mathcal{N})}
        \arrow[r,"(\varphi^*)_*"']
    & {\mathcal{Z}_1(\mathcal{M}) \times_{\mathcal{Z}_1(G;H)} \mathcal{Z}_1(\mathcal{N})}
\end{tikzcd},\] where:
\begin{itemize}
    \item $\mathcal{Z}_1(G;H)$ is the monoidal intertwiner of $G$ and $H$ (defined below);
    
    \item The right bottom corner is the 2-fiber product of \emph{categories}: \[\begin{tikzcd}
        {\mathcal{Z}_1(\mathcal{M})}
            \arrow[dd,"G_*"']
        & {\mathcal{Z}_1(\mathcal{M}) \times_{\mathcal{Z}_1(G;H)} \mathcal{Z}_1(\mathcal{N})}
            \arrow[r,"\mathbf{pr2}"]
            \arrow[l,"\mathbf{pr1}"']
        & {\mathcal{Z}_1(\mathcal{N})}
            \arrow[dd,"H^*"]
        \\ {}
            \arrow[rr,Rightarrow,shorten <=75pt, shorten >=75pt]
        & {}
        & {}
        \\ {\mathcal{Z}_1(G)}
            \arrow[r,"\varphi_*"']
        & {\mathcal{Z}_1(G;H)}
        & {\mathcal{Z}_1(H)}
            \arrow[l,"\varphi^*"]
    \end{tikzcd};\]
    
    \item $\varphi_* \colon \mathcal{Z}_1(G) \to \mathcal{Z}_1(G;H); \quad (x,\beta^x) \mapsto (x,(\varphi \otimes x) \circ_2 \beta^x);$
    
    \item $\varphi^* \colon \mathcal{Z}_1(H) \to \mathcal{Z}_1(G;H); \quad (y,\beta^y) \mapsto (y,\beta^y \circ_2 (y \otimes \varphi)).$
\end{itemize}

\begin{Definition}
    Given two monoidal functors $G,H \colon \mathcal{M} \to \mathcal{N}$, the \emph{monoidal intertwiner} $\mathcal{Z}_1(G;H)$ is the category of \emph{lax} half-braidings in $\mathcal{N}$ intertwining $G$ and $H$:

        \begin{enumerate}
                \item An object is a pair $(x,\beta^x)$, where \begin{enumerate}
                    \item $x$ is an object in $\mathcal{N}$,
                    
                    \item $\beta^x_y:x \otimes_\mathcal{N} G(y) \to H(y) \otimes_\mathcal{N} x$ is a natural transformation given on $y$ in $\mathcal{M}$, \end{enumerate}
                    such that the following hexagon commutes naturally in $y,z$: \[\begin{tikzcd}[column sep=40pt]
                    {x \otimes_\mathcal{N} G(y \otimes_\mathcal{M} z)}
                        \arrow[r,"\beta^x_{y \otimes z}"]
                    & {H(y \otimes_\mathcal{M} z) \otimes_\mathcal{N} x}
                        \arrow[d,"\gamma^H_{y,z} \otimes x"]
                    \\ {x \otimes_\mathcal{N} (G(y) \otimes_\mathcal{N} G(z))}
                        \arrow[u,"x \otimes (\gamma^G_{y,z})^{-1}"]
                    & {(H(y) \otimes_\mathcal{N} H(z)) \otimes_\mathcal{N} x}
                        \arrow[d,"\alpha_{Hy,Hz,x}"]
                    \\ {(x \otimes_\mathcal{N} G(y)) \otimes_\mathcal{N} G(z)}
                        \arrow[d,"\beta^x_{y} \otimes G(z)"']
                        \arrow[u,"\alpha_{x,Gy,Gz}"]
                    & {H(y) \otimes_\mathcal{N} (H(z) \otimes_\mathcal{N} x)}
                    \\ {(H(y) \otimes_\mathcal{N} x) \otimes_\mathcal{N} G(z)}
                        \arrow[r,"\alpha_{Hy,x,Gz}"']
                    & {H(y) \otimes_\mathcal{N} (x \otimes_\mathcal{N} G(z))}
                        \arrow[u,"H(y) \otimes \beta^x_{z}"']
                    \end{tikzcd}.\]

                \item A morphism from $(x,\beta^x)$ to $(y,\beta^y)$ consists of a morphism $f:x \to y$ in $\mathcal{N}$ such that the diagram commutes naturally in $z$: \[\begin{tikzcd}[sep=large]
                {x \otimes_\mathcal{N} G(z)}
                    \arrow[r,"\beta^{x}_{z}"]
                    \arrow[d,"f \otimes G(z)"']
                & {H(z) \otimes_\mathcal{N} x}
                \arrow[d,"H(z) \otimes f"]
                \\ {y \otimes_\mathcal{N} G(z)}
                    \arrow[r,"\beta^{y}_{z}"']
                & {H(z) \otimes_\mathcal{N} y}
                \end{tikzcd}. \]
            \end{enumerate}
\end{Definition}

\begin{Remark}
    We emphasize that here we must allow lax half-braidings for the construction of $\varphi_*$ and $\varphi^*$, since the natural transformation $\varphi$ need not be invertible. In particular, $\mathcal{Z}_1(G,G)$ is generally larger than the monoidal centralizer $\mathcal{Z}_1(G)$.
\end{Remark}

\begin{Remark}
    We can characterize the full subcategory of $\mathcal{Z}_1(G;H)$ consisting of objects with invertible half-braidings as \[\mathbf{Fun}_{\mathcal{M} \boxtimes \mathcal{N}^{mp}}({}_{\langle G \rangle} \mathcal{N},{}_{\langle H \rangle} \mathcal{N}) \simeq \mathbf{Fun}_{\mathcal{M} \boxtimes \mathcal{M}^{mp}}(\mathcal{M},{}_{\langle H \rangle} \mathcal{N}_{\langle G \rangle}).\]
\end{Remark}

\begin{Remark}
    The monoidal intertwiner $\mathcal{Z}_1(G;H)$ is not necessarily a monoidal category, but it has a right action of $\mathcal{Z}_1(G)$ and a left action of $\mathcal{Z}_1(H)$.
\end{Remark}

\begin{Lemma}
    $\mathcal{Z}_1(\mathcal{M}) \times_{\mathcal{Z}_1(G;H)} \mathcal{Z}_1(\mathcal{N})$ inherits a canonical monoidal structure from $\mathcal{Z}_1(G)$ and $\mathcal{Z}_1(H)$. Meanwhile, $\mathcal{Z}_1(\mathcal{M}) \times^\varphi \mathcal{Z}_1(\mathcal{N})$ inherits a canonical braiding from $\mathcal{Z}_1(\mathcal{M})$ and $\mathcal{Z}_1(\mathcal{N})$.
\end{Lemma}

\begin{proof}
    Apply Theorem \ref{thm:AlgebraStructureOn2FiberProduct}.
\end{proof}

\begin{Remark}
    A generic object in $\mathcal{Z}_1(\mathcal{M}) \times^\varphi \mathcal{Z}_1(\mathcal{N})$ consists of \begin{enumerate}
        \item A half-braiding $(m,\beta^m)$ in $\mathcal{M}$, 
        
        \item A half-braiding $(n,\beta^n)$ in $\mathcal{N}$,
        
        \item An invertible 2-morphism $\xi^G:G(m) \to n$ in $\mathcal{N}$,
        
        \item An invertible 2-morphism $\xi^H:H(m) \to n$ in $\mathcal{N}$,
    \end{enumerate} satisfying the following conditions:
    \begin{enumerate}
        \item[a.] $(m,\beta^m)$, $(n,\beta^n)$ and $\xi^G$ satisfy the coherence condition in Remark \ref{rmk:CoherenceForMonoidalCentralizer} for the monoidal functor $G$;
        
        \item[b.] $(m,\beta^m)$, $(n,\beta^n)$ and $\xi^H$ satisfy the coherence condition in Remark \ref{rmk:CoherenceForMonoidalCentralizer} for the monoidal functor $H$;
        
        \item[c.] $\xi^H \circ_1 \varphi_m = \xi^G$, i.e., the coherence condition (9) in \cite[Definition 3.3]{HPT23}.
    \end{enumerate}
\end{Remark}

\begin{Proposition}
    A monoidal natural transformation $\varphi \colon G \to H$ of central $\mathcal{A}$-module monoidal functors $(G,\psi^G)$ and $(H,\psi^H)$ preserves central $\mathcal{A}$-module monoidal structures if and only if there exists a braided monoidal functor \[F^\varphi:\mathcal{A} \to \mathcal{Z}_1(\mathcal{M}) \times^\varphi \mathcal{Z}_1(\mathcal{N})\]
    \[x \mapsto \left( F^\mathcal{M}(x), F^\mathcal{N}(x), G(F^\mathcal{M}(x)) \xrightarrow{\psi^G_x} F^\mathcal{N}(x), H(F^\mathcal{M}(x)) \xrightarrow{\psi^H_x} F^\mathcal{N}(x) \right).\]
\end{Proposition}

\begin{proof}
    Unpacking the definition of a central $\mathcal{A}$-module monoidal natural transformation, the additional coherence condition requires the following two natural isomorphisms to coincide:
        \[\begin{tikzcd}
            {\mathcal{Z}_1(\mathcal{M})}
                \arrow[dd,"G_*"']
                \arrow[rd,Rightarrow,shorten <=10pt, shorten >=10pt,"\phi^G"']
            & {\mathcal{A}}
                \arrow[r,"F^\mathcal{N}"]
                \arrow[l,"F^\mathcal{M}"']
            & {\mathcal{Z}_1(\mathcal{N})}
                \arrow[dd,"H^*"]
                \arrow[ddll,"G^*"]
            \\ {}
            & {}
                \arrow[rd,Rightarrow,shorten <=20pt, shorten >=10pt,"\mathrm{id}"]
            & {}
            \\ {\mathcal{Z}_1(G)}
                \arrow[r,"\varphi_*"']
            & {\mathcal{Z}_1(G;H)}
            & {\mathcal{Z}_1(H)}
                \arrow[l,"\varphi^*"]
        \end{tikzcd}\] and\[\begin{tikzcd}
            {\mathcal{Z}_1(\mathcal{M})}
                \arrow[dd,"G_*"']
                \arrow[ddrr,"H_*"']
            & {\mathcal{A}}
                \arrow[r,"F^\mathcal{N}"]
                \arrow[l,"F^\mathcal{M}"']
            & {\mathcal{Z}_1(\mathcal{N})}
                \arrow[dd,"H^*"]
            \\ {}
            & {}
                \arrow[ur,Rightarrow,shorten <=10pt, shorten >=20pt,"\phi^H"']
            & {}
            \\ {\mathcal{Z}_1(G)}
                \arrow[r,"\varphi_*"']
                \arrow[ru,Rightarrow,shorten <=10pt, shorten >=20pt,"\mathrm{id}"]
            & {\mathcal{Z}_1(G;H)}
            & {\mathcal{Z}_1(H)}
                \arrow[l,"\varphi^*"]
        \end{tikzcd}.\]

    By the 2-universal property of the 2-fiber product of categories, this induces a functor $F^\varphi \colon \mathcal{A} \to \mathcal{Z}_1(\mathcal{M}) \times^\varphi \mathcal{Z}_1(\mathcal{N})$. Finally, this canonical functor preserves the monoidal structures and braidings in the same manner as described in Theorem \ref{thm:AlgebraStructureOn2FiberProduct}.
\end{proof}

In summary, a central $\mathcal{A}$-module monoidal natural transformation 
\[\begin{tikzcd}[column sep=40pt]
        {\mathcal{M}}
            \arrow[r,bend left=50pt,"{(G,\psi^G)}"{name=U}]
            \arrow[r,bend right=50pt,"{(H,\psi^H)}"'{name=D}]
        & {\mathcal{N}}
        \arrow[Rightarrow, from=U, to=D, "\varphi", shorten <=10pt, shorten >=10pt]
\end{tikzcd} \] corresponds to a 2-span of braided monoidal categories under $\mathcal{A}$:
\[\begin{tikzcd}
    {}
    & {}
    & {\mathcal{Z}_1(\mathcal{M}) \times_{\mathcal{Z}_1(G)} \mathcal{Z}_1(\mathcal{N})}
        \arrow[lldd]
        \arrow[rrdd]
    & {}
    & {}
    \\ {}
    & {}
    & {}
    & {}
    & {}
    \\ {\mathcal{Z}_1(\mathcal{M})}
    & {}
    & {\mathcal{Z}_1(\mathcal{M}) \times^\varphi \mathcal{Z}_1(\mathcal{N})}
        \arrow[uu]
        \arrow[dd]
    & {}
    & {\mathcal{Z}_1(\mathcal{N})}
    \\ {}
    & {}
    & {}
    & {}
    & {}
    \\ {}
    & {}
    & {\mathcal{Z}_1(\mathcal{M}) \times_{\mathcal{Z}_1(H)} \mathcal{Z}_1(\mathcal{N})}
        \arrow[lluu]
        \arrow[rruu]
    & {}
    & {}
\end{tikzcd}.\]

\begin{Corollary}
    The above construction provides a lax 3-functor from the 2-category of central $\mathcal{A}$-module monoidal categories to the 3-category of 2-spans of braided monoidal categories under $\mathcal{A}$.
\end{Corollary}

\vspace*{1cm}

Let us consider a similar example when $\mathfrak{C} = \mathbf{BrCat}$ is the 2-category of braided monoidal categories with the Cartesian product. Internal homs do not exist in $\mathfrak{C}$.

\begin{Remark}
    Algebras in $\mathfrak{C}$ are symmetric monoidal categories, or equivalently categories equipped with three compatible monoidal structures \cite{BFSV03}, a fact which further generalizes the Eckmann-Hilton arguments. Notice that symmetric monoidal categories are \emph{stable} in the sense that adding further compatible monoidal structures does not provide any new information, see \emph{Stabilization Hypothesis} in \cite{BD}.
\end{Remark}

\medskip \noindent \textit{On the level of Objects.} Given a symmetric monoidal category $\mathcal{S}$, a left module of $\mathcal{S}$ in $\mathfrak{C}$ is a central $\mathcal{S}$-module braided monoidal category.

\begin{Definition}
    An $\mathcal{S}$-central braided monoidal category consists of a braided category $\mathcal{A}$ and a braided functor $F^\mathcal{A}:\mathcal{S} \to \mathcal{Z}_2(\mathcal{A})$, where $\mathcal{Z}_2(\mathcal{A})$ is the M{\"u}ger center of $\mathcal{A}$.
\end{Definition}

\begin{Remark}
    The \emph{M{\"u}ger center} \cite{Mu3} or the \emph{symmetric center} $\mathcal{Z}_2(\mathcal{A})$ is the full subcategory of $\mathcal{A}$ consisting of objects $x$ such that for any object $y$ in $\mathcal{A}$, the double braiding $\beta^{y}_x \circ \beta^{x}_y$ is the identity morphism on $x \otimes y$. It is a symmetric monoidal category.
\end{Remark}

\begin{Definition}
    By analogy with the monoidal centralizer, for a braided functor $G \colon \mathcal{A} \to \mathcal{B}$, we can define the \emph{braided centralizer} $\mathcal{Z}_2(G)$ as the full subcategory of $\mathcal{B}$ consisting of objects $x$ whose double braiding with any object in the image of $G$ is trivial. Note that $\mathcal{Z}_2(G)$ is only a braided monoidal category in general.

    For two braided functors $G,H \colon \mathcal{A} \to \mathcal{B}$, we can similarly define the \emph{braided intertwiner} $\mathcal{Z}_2(G;H)$ as the smallest full subcategory of $\mathcal{B}$ containing $\mathcal{Z}_2(G)$ and $\mathcal{Z}_2(H)$.    
\end{Definition}

\medskip \noindent \textit{On the level of 1-Morphisms.} Given two central $\mathcal{S}$-module braided categories $\mathcal{A}$ and $\mathcal{B}$, a braided functor $G \colon \mathcal{A} \to \mathcal{B}$ preserves the central $\mathcal{S}$-module braided monoidal structures if it is equipped with a monoidal natural isomorphism \[\begin{tikzcd}[column sep=large]
    {}
    & {\mathcal{Z}_2(\mathcal{A})}
        \arrow[r,hook]
        \arrow[dd,Rightarrow,shorten <=10pt,shorten >=10pt,"\psi^G"]
    & {\mathcal{A}}
        \arrow[dd,"G"]
    \\  {\mathcal{S}}
        \arrow[ur,"F^\mathcal{A}"]
        \arrow[dr,"F^\mathcal{B}"']
    & {}
    & {}
    \\ {}
    & {\mathcal{Z}_2(\mathcal{B})}
        \arrow[r,hook]
    & {\mathcal{B}}
\end{tikzcd}.\]

\begin{Proposition}
    The data of a central $\mathcal{S}$-module braided functor structure $\psi^G$ corresponds to a braided functor \[\mathcal{S} \to \mathcal{Z}_2(\mathcal{A}) \times_{\mathcal{Z}_2(G)} \mathcal{Z}_2(\mathcal{B}); \quad x \mapsto \left( F^\mathcal{A}(x), F^\mathcal{B}(x), G(F^\mathcal{A}(x)) \xrightarrow{\psi^G_x} F^\mathcal{B}(x) \right),\] where the 2-fiber product is defined via
    \[\begin{tikzcd}
    {\mathcal{Z}_2(\mathcal{A}) \times_{\mathcal{Z}_2(G)} \mathcal{Z}_2(\mathcal{B})}
        \arrow[r,"\mathbf{pr1}"]
        \arrow[d,"\mathbf{pr2}"']
    & {\mathcal{Z}_2(\mathcal{A})}
        \arrow[d,"G_*"]
        \arrow[dl,Rightarrow,shorten <=25pt, shorten >=20pt]
    \\ {\mathcal{Z}_2(\mathcal{B})}
        \arrow[r,"G^*"']
    & {\mathcal{Z}_2(G)}
\end{tikzcd}.\]
\end{Proposition}

\begin{proof}
    Observe that $\psi^G$ can be rewritten as a monoidal natural isomorphism \[\begin{tikzcd}[sep=large]
    {\mathcal{S}}
        \arrow[r,"F^\mathcal{A}"]
        \arrow[d,"F^\mathcal{B}"']
    & {\mathcal{Z}_2(\mathcal{A})}
        \arrow[d,"G_*"]
        \arrow[dl,Rightarrow,shorten <=20pt, shorten >=15pt,"\phi^G"']
    \\ {\mathcal{Z}_2(\mathcal{B})}
        \arrow[r,"G^*"']
    & {\mathcal{Z}_2(G)}
\end{tikzcd}.\] Then the rest follows from the 2-universal property of the 2-fiber product.
\end{proof}

\begin{Remark}
    The 2-fiber product $\mathcal{Z}_2(\mathcal{A}) \times_{\mathcal{Z}_2(G)} \mathcal{Z}_2(\mathcal{B})$ is a full subcategory of $\mathcal{Z}_2(\mathcal{A})$, hence it is always symmetric. 
\end{Remark}

\medskip \noindent \textit{On the level of 2-Morphisms.} A central $\mathcal{S}$-module monoidal natural transform between $(G,\psi^G)$ and $(H,\psi^H)$ is simply a monoidal natural transform $\varphi:G \to H$ such that
\[\begin{tikzcd}
    {}
    & {\mathcal{Z}_2(\mathcal{A})}
        \arrow[r,hook]
        \arrow[dd,Rightarrow,shorten <=10pt,shorten >=10pt,"\psi^G"]
    & {\mathcal{A}}
        \arrow[dd,bend left=35pt,"H"]
        \arrow[dd,bend right=35pt,"G"']
    \\  {\mathcal{S}}
        \arrow[ur,"F^\mathcal{A}"]
        \arrow[dr,"F^\mathcal{B}"']
    & {}
    & {}
        \arrow[l,Leftarrow,xshift=25pt,shorten <=15pt,shorten >=15pt,"\varphi"']
    \\ {}
    & {\mathcal{Z}_2(\mathcal{B})}
        \arrow[r,hook]
    & {\mathcal{B}}
\end{tikzcd} \quad = \quad \begin{tikzcd}
    {}
    & {\mathcal{Z}_2(\mathcal{A})}
        \arrow[r,hook]
        \arrow[dd,Rightarrow,shorten <=10pt,shorten >=10pt,"\psi^H"]
    & {\mathcal{A}}
        \arrow[dd,"H"]
    \\  {\mathcal{S}}
        \arrow[ur,"F^\mathcal{A}"]
        \arrow[dr,"F^\mathcal{B}"']
    & {}
    & {}
    \\ {}
    & {\mathcal{Z}_2(\mathcal{B})}
        \arrow[r,hook]
    & {\mathcal{B}}
\end{tikzcd}.\]

Consider the 2-fiber product
\[\begin{tikzcd}
    {\mathcal{Z}_2(\mathcal{A}) \times^\varphi \mathcal{Z}_2(\mathcal{B})}
        \arrow[r,"\mathbf{pr1}"]
        \arrow[d,"\mathbf{pr2}"']
    & {\mathcal{Z}_2(\mathcal{A}) \times_{\mathcal{Z}_2(G)} \mathcal{Z}_2(\mathcal{B})}
        \arrow[d,"(\varphi_*)^*"]
        \arrow[dl,Rightarrow,shorten <=25pt, shorten >=20pt]
    \\ {\mathcal{Z}_2(\mathcal{A}) \times_{\mathcal{Z}_2(H)} \mathcal{Z}_2(\mathcal{B})}
        \arrow[r,"(\varphi_*)^*"']
    & {\mathcal{Z}_2(\mathcal{A}) \times_{\mathcal{Z}_2(G;H)} \mathcal{Z}_2(\mathcal{B})}
\end{tikzcd},\] where:

\begin{itemize}
    \item $\varphi_* \colon \mathcal{Z}_2(G) \to \mathcal{Z}_2(G;H)$ and $\varphi^* \colon \mathcal{Z}_2(H) \to \mathcal{Z}_2(G;H)$ are canonically defined embeddings;
    
    \item The right bottom corner is the 2-fiber product \[\begin{tikzcd}
        {\mathcal{Z}_2(\mathcal{A})}
            \arrow[dd,"G_*"']
        & {\mathcal{Z}_2(\mathcal{A}) \times_{\mathcal{Z}_2(G;H)} \mathcal{Z}_2(\mathcal{B})}
            \arrow[r,"\mathbf{pr2}"]
            \arrow[l,"\mathbf{pr1}"']
        & {\mathcal{Z}_2(\mathcal{B})}
            \arrow[dd,"H^*"]
        \\ {}
            \arrow[rr,Rightarrow,shorten <=75pt, shorten >=75pt]
        & {}
        & {}
        \\ {\mathcal{Z}_2(G)}
            \arrow[r,"\varphi_*"']
        & {\mathcal{Z}_2(G;H)}
        & {\mathcal{Z}_2(H)}
            \arrow[l,"\varphi^*"]
    \end{tikzcd}.\]
\end{itemize}

\begin{Remark}
    A generic object in $\mathcal{Z}_2(\mathcal{A}) \times^\varphi \mathcal{Z}_2(\mathcal{B})$ consists of a transparent object $x$ in $\mathcal{A}$ and a transparent object $y$ in $\mathcal{B}$, together with an invertible 2-morphisms $\xi^G:G(x) \to y$ and $\xi^H:H(x) \to y$ in $\mathcal{B}$, such that $\xi^H \circ_1 \varphi_x = \xi^G$.
    
    A morphism from $(x_0,y_0,\xi^G_0,\xi^H_0) \to (x_1,y_1,\xi^G_1,\xi^H_1)$ consists of a morphism $f:x_0 \to x_1$ in $\mathcal{A}$ and a morphism $g:y_0 \to y_1$ in $\mathcal{B}$, such that the following diagram commutes:
    \[\begin{tikzcd}[sep=large]
        {G(x_0)}
            \arrow[r,"G(f)"]
            \arrow[d,"\xi^G_0"']
        & {G(x_1)}
            \arrow[d,"\xi^G_1"]
        \\ {y_0}
            \arrow[r,"g"']
        & {y_1}
    \end{tikzcd}, \quad \begin{tikzcd}[sep=large]
        {H(x_0)}
            \arrow[r,"H(f)"]
            \arrow[d,"\xi^H_0"']
        & {H(x_1)}
            \arrow[d,"\xi^H_1"]
        \\ {y_0}
            \arrow[r,"g"']
        & {y_1} 
    \end{tikzcd}.\]
\end{Remark}

\begin{Remark}
    $\mathcal{Z}_2(\mathcal{A}) \times^\varphi \mathcal{Z}_2(\mathcal{B})$ is the largest full subcategory of $\mathcal{B}$ contained in both $\mathcal{Z}_2(G)$ and $\mathcal{Z}_2(H)$.
\end{Remark}

\begin{Proposition}
    A monoidal natural transformation $\varphi \colon G \to H$ preserves central $\mathcal{S}$-module braided monoidal structures if and only if there exists a braided monoidal functor \[\mathcal{S} \to \mathcal{Z}_2(\mathcal{A}) \times^\varphi \mathcal{Z}_2(\mathcal{B}); \] 
    \[x \mapsto \left( F^\mathcal{A}(x), F^\mathcal{B}(x), G(F^\mathcal{A}(x)) \xrightarrow{\psi^G_x} F^\mathcal{B}(x), H(F^\mathcal{A}(x)) \xrightarrow{\psi^H_x} F^\mathcal{B}(x) \right).\]
\end{Proposition}

\begin{proof}
    Observe the following equal monoidal natural isomorphisms: \[\begin{tikzcd}
        {\mathcal{Z}_2(\mathcal{A})}
            \arrow[dd,"G_*"']
            \arrow[rd,Rightarrow,shorten <=10pt, shorten >=10pt,"\phi^G"']
        & {\mathcal{S}}
            \arrow[r,"F^\mathcal{B}"]
            \arrow[l,"F^\mathcal{A}"']
        & {\mathcal{Z}_2(\mathcal{B})}
            \arrow[dd,"H^*"]
            \arrow[ddll,"G^*"]
        \\ {}
        & {}
            \arrow[rd,Rightarrow,shorten <=15pt, shorten >=10pt,"\mathrm{id}"]
        & {}
        \\ {\mathcal{Z}_2(G)}
            \arrow[r,"\varphi_*"']
        & {\mathcal{Z}_2(G;H)}
        & {\mathcal{Z}_2(H)}
            \arrow[l,"\varphi^*"]
        \end{tikzcd}\] and \[\begin{tikzcd}
        {\mathcal{Z}_2(\mathcal{A})}
            \arrow[dd,"G_*"']
            \arrow[ddrr,"H_*"']
        & {\mathcal{S}}
            \arrow[r,"F^\mathcal{B}"]
            \arrow[l,"F^\mathcal{A}"']
        & {\mathcal{Z}_2(\mathcal{B})}
            \arrow[dd,"H^*"]
        \\ {}
        & {}
            \arrow[ur,Rightarrow,shorten <=10pt, shorten >=15pt,"\phi^H"']
        & {}
        \\ {\mathcal{Z}_2(G)}
            \arrow[r,"\varphi_*"']
            \arrow[ru,Rightarrow,shorten <=10pt, shorten >=20pt,"\mathrm{id}"]
        & {\mathcal{Z}_2(G;H)}
        & {\mathcal{Z}_2(H)}
            \arrow[l,"\varphi^*"]
        \end{tikzcd}.\]
        The result then follows from the 2‑universal property of the 2‑fiber product.
\end{proof}

In summary, a central $\mathcal{S}$-module monoidal natural transformation 
\[\begin{tikzcd}[column sep=40pt]
    {\mathcal{A}}
        \arrow[r,bend left=50pt,"{(G,\psi^G)}"{name=U}]
        \arrow[r,bend right=50pt,"{(H,\psi^H)}"'{name=D}]
    & {\mathcal{B}}
    \arrow[Rightarrow, from=U, to=D, "\varphi", shorten <=10pt, shorten >=10pt]
\end{tikzcd} \] corresponds to a 2-span of symmetric monoidal categories under $\mathcal{S}$:
\[\begin{tikzcd}
    {}
    & {}
    & {\mathcal{Z}_2(\mathcal{A}) \times_{\mathcal{Z}_2(G)} \mathcal{Z}_2(\mathcal{B})}
    \arrow[lldd]
    \arrow[rrdd]
    & {}
    & {}
    \\ {}
    & {}
    & {}
    & {}
    & {}
    \\ {\mathcal{Z}_2(\mathcal{A})}
    & {}
    & {\mathcal{Z}_2(\mathcal{A}) \times^\varphi \mathcal{Z}_2(\mathcal{B})}
    \arrow[uu]
    \arrow[dd]
    & {}
    & {\mathcal{Z}_2(\mathcal{B})}
    \\ {}
    & {}
    & {}
    & {}
    & {}
    \\ {}
    & {}
    & {\mathcal{Z}_2(\mathcal{A}) \times_{\mathcal{Z}_2(H)} \mathcal{Z}_2(\mathcal{B})}
    \arrow[lluu]
    \arrow[rruu]
    & {}
    & {}
\end{tikzcd}.\]

\begin{Corollary}
    The above construction provides a lax 3-functor from the 2-category of central $\mathcal{S}$-module braided monoidal categories to the 3-category of 2-spans of symmetric monoidal categories under $\mathcal{S}$.
\end{Corollary}